\documentclass[11pt,reqno]{amsart}
\usepackage{euscript,amssymb}
\usepackage[arrow,matrix,curve]{xy}
\usepackage{array,longtable}
\usepackage{a4wide}

\newenvironment{m-theorem}{%
\vskip5pt\refstepcounter{stff}\trivlist \itemindent 0pt
\item[\hskip\labelsep\bf Theorem~\thestff]%
\it\ignorespaces}{\endtrivlist\vskip5pt}%
\newenvironment{m-proposition}{%
\vskip5pt\refstepcounter{stff}\trivlist \itemindent 0pt
\item[\hskip\labelsep\bf Proposition~\thestff]%
\it\ignorespaces}{\endtrivlist\vskip5pt}%
\newenvironment{m-corollary}{%
\vskip5pt\refstepcounter{stff}\trivlist \itemindent 0pt
\item[\hskip\labelsep\bf Corollary~\thestff]%
\it\ignorespaces}{\endtrivlist\vskip5pt}%
\newenvironment{m-lemma}{%
\vskip5pt\refstepcounter{stff}\trivlist \itemindent 0pt
\item[\hskip\labelsep\bf Lemma~\thestff]%
\it\ignorespaces}{\endtrivlist\vskip5pt}%
\newenvironment{m-definition}{%
\vskip5pt\refstepcounter{stff}\trivlist \itemindent 0pt
\item[\hskip\labelsep\bf Definition~\thestff]%
\ignorespaces}{\endtrivlist\vskip5pt}%
\newenvironment{m-notation}{%
\vskip5pt\refstepcounter{stff}\trivlist \itemindent 0pt
\item[\hskip\labelsep\bf Notation~\thestff]%
\ignorespaces}{\endtrivlist\vskip5pt}%
\newenvironment{m-example}{%
\vskip5pt\refstepcounter{stff}\trivlist \itemindent 0pt
\item[\hskip\labelsep\bf Example~\thestff]%
\ignorespaces}{\endtrivlist\vskip5pt}
\newenvironment{m-remark}{%
\vskip5pt\refstepcounter{stff}\trivlist \itemindent 0pt
\item[\hskip\labelsep\bf Remark~\thestff]%
\ignorespaces}{\endtrivlist\vskip5pt}
\newenvironment{m-question}{%
\vskip5pt\refstepcounter{stff}\trivlist \itemindent 0pt
\item[\hskip\labelsep\bf Question.]%
\ignorespaces}{\endtrivlist\vskip5pt}%
\newenvironment{thm-nono}{
\vskip5pt\trivlist \itemindent 0pt
\item[\hskip\labelsep\bf Theorem.]%
\it\ignorespaces}{\endtrivlist\vskip5pt}%
\newenvironment{m-thank}{%
\vskip5pt\trivlist \itemindent 0pt
\item[\hskip\labelsep\it Acknowledgments]%
\ignorespaces}{\endtrivlist\vskip5pt}%

\let\hra\hookrightarrow
\let\mt\mapsto

\font\tenmsa=msam10 %
\newcommand\hdashpiece{%
{\vrule height2.75pt depth-2.35pt width2.3pt \kern1.7pt}}%
\newcommand\hdashpieces{%
{\hdashpiece\hdashpiece\hdashpiece\hdashpiece}}%
\let\dashto\dashrightarrow
\newcommand\dashar{\mathrel{%
\hdashpieces\kern-0.4pt\hbox{\tenmsa K}}}%

\let\euf\EuScript 
\let\cal\mathcal
\let\mbb\mathbb

\DeclareFontFamily{OT1}{rsfs}{}
\DeclareFontShape{OT1}{rsfs}{n}{it}{<->rsfs10}{}
\DeclareMathAlphabet{\crl}{OT1}{rsfs}{n}{it}

\newcommand\ouset[3]{{\overset{#2}{\underset{#1}#3}}}
\let\ovl\overline
\let\unl\underline
\let\unbar\underbar

\let\tld\tilde
\let\wtld\widetilde

\let\nit\noindent

\let\disp\displaystyle
\let\srel\stackrel
\let\vphi\varphi

\let\veps\varepsilon

\let\lan\langle
\let\ran\rangle
\newcommand\lran[1]{{\lan #1\ran}}

\newcommand\rd{{\rm d}} 
\newcommand\ev{{\rm ev}}
\newcommand\Aut{\operatorname{\textrm{Aut}\kern1pt}}
\newcommand\cAut{\operatorname{\mathcal{A}\kern-1pt\textit{ut}\kern1pt}}
\newcommand\End{\operatorname{\rm{End}\kern1pt}}
\newcommand\cEnd{\operatorname{\mathcal{E}\kern-1pt\textit{nd}\kern1pt}}
\newcommand\Ext{\mathop{\rm Ext}\nolimits}
\newcommand\Flag{{\rm Flag}}

\newcommand\Hom{\mathop{\rm Hom}\nolimits}
\newcommand\cHom{\operatorname{\mathcal{H}\kern-1pt\textit{om}\kern1pt}}
\newcommand\Img{{\rm Im}}
\newcommand\Ker{{\rm Ker}}

\newcommand\Lie{\mathop{\rm Lie}\nolimits}

\newcommand{\Weyl}{\mathop{\rm{Weyl}}\nolimits}
\newcommand\NS{\mathop{\rm NS}\nolimits}
\newcommand\Pic{\mathop{\rm Pic}\nolimits}

\newcommand\Spec{\mathop{\rm Spec}\nolimits}

\newcommand\Sym{\mathop{\rm Sym}\nolimits}
\newcommand\invq{{\slash\kern-2.5pt\slash}}
\newcommand\pr{\mathop{\rm pr}\nolimits}
\newcommand\rk{{\rm rk}}

\newcommand\Grs{{\rm Gr}}
\newcommand\spGrs{{\rm sp{\hdashpiece\!}Gr}}
\newcommand\oGrs{{\rm o{\hdashpiece\!}Gr}}

\numberwithin{equation}{section}
\numberwithin{figure}{section} 

\let\l\lambda
\let\L\Lambda

\let\si\sigma

\let\sm\setminus

\newcommand\bbC{{\mbb C}}

\newcommand\bbF{{\mbb F}}

\newcommand\bbk{\mbox{\rm I\kern-1.5pt k}}
\newcommand\sbbk{\hbox{\scriptsize I{\kern-.8pt}k}}

\newcommand\bbQ{{\mbb Q}}
\newcommand\bbR{{\mbb R}}

\newcommand\bbY{{\mbb Y}}
\newcommand\bone{{1\kern-0.57ex\rm l}}

\newcommand\bk{{\mbb C}}

\newcommand\eA{{\euf A}}
\newcommand\eB{{\euf B}}
\newcommand\eC{{\euf C}}
\newcommand\cE{{\cal E}}

\newcommand\eF{{\euf F}}

\newcommand\eI{{\euf I}}
\newcommand\eJ{{\euf J}}

\newcommand\eL{{\euf L}}

\newcommand\eN{{\euf N}}

\newcommand\eO{{\euf O}}

\newcommand\cU{{\cal U}}
\newcommand\eU{{\euf U}}
\newcommand\cY{{\cal Y}}

\newcommand\eT{{\euf T}}
\newcommand\eW{{\euf W}}

\newcommand\uz{{\unl z}}

\newcommand\codim{{\rm codim}}

\newcommand\Ad{{\rm Ad}}
\newcommand\Bl{{\rm Bl}}
\newcommand\GL{{\rm GL}}

\newcommand\SO{{\rm SO}}

\newcommand\Sp{{\rm Sp}}

\renewcommand\det{{\rm det}}

\let\ges\geqslant
\let\les\leqslant

\newcommand\Supp{{\rm Supp}\,}

\newcommand\res{{\rm res}}
\let\surj\twoheadrightarrow

\newcommand\pos{{>0}}
\newcommand\apos{{\,\gtrsim0}}
\newcommand\ball {{\mbb B}}
\newcommand\BB{{\rm BB}}
\newcommand\Bht{{\rm Bruhat}}
\newcommand\Levi{\mathop{\rm Levi}\nolimits}
\newcommand\diag{\mathop{\rm diag}\nolimits}
\newcommand{\notr}{\mbox{(no-$\triangle$)\ }}
\newcommand{\arm}{\mbox{($1$-arm)\ }}
\newcommand\cd{{\rm cd}}
\newcommand\ysx{{\rm(YSX)}{\hskip1ex}}

\author{Mihai Halic}
\address{}


\begin{document}

\title{Vector bundles on projective varieties\\ which split along $q$-ample subvarieties}

\begin{abstract}
Given a vector bundle on a complex, smooth projective variety, I prove that 
its splitting along a very general, $q$-ample subvariety (for appropriate $q$), 
which admits sufficiently many embedded deformations implies the global splitting. 
The result goes significantly beyond the previously known splitting criteria 
obtained by restricting vector bundles to subvarieties. 

I discuss the particular case of zero loci of sections in globally generated vector bundles, 
on one hand, and of sources of multiplicative group actions 
(corresponding to Bialynicki-Birula decompositions), on the other hand. 
Finally, I elaborate on the symplectic and orthogonal Grassmannians; 
I prove that the splitting of an arbitrary vector bundle on them can be read off 
from its restriction to a low dimensional `sub'-Grassmannian.
\end{abstract}

\maketitle
\markboth{\sc Mihai Halic}%
{\sc Vector bundles which split along $q$-ample subvarieties}

\section*{Introduction}

For  a vector bundle $\crl V$ on the irreducible projective variety $X$, I consider its 
restriction $\crl V_Y$ to an irreducible subvariety $Y\subset X$, and investigate the 
following: 
\begin{m-question}
Assuming that $\crl V_Y$ splits, under what assumptions on $Y$ does $\crl V$ split too?
\end{m-question}
To my knowledge, 
\emph{there is not a single reference which addresses the question in this generality}. 
Yet, the problem is interesting because it allows to probe the splitting of vector bundles on 
(high dimensional) varieties by restricting them to (possibly low dimensional) subvarieties. 
The goal of this article is to give a tentative answer to this question. 

Horrocks' criterion \cite{ho} is the most widely known result of this type. 
It was extended in \cite{ba}, for restrictions to ample divisors on (suitable) varieties. 
Surprisingly, restrictions to higher codimensional subvarieties have not been considered 
yet. For this reason, I studied in \cite{hal} the case when $Y$ is the zero locus of 
a regular section in a globally generated, \emph{ample} vector bundle $\eN$ on $X$. 

Unfortunately, this setting discards very basic situations, \textit{e.g.} $X$ is a product 
$X'\times V$ and $\eN=\eO_{X'}(1)$. 
Second, in perfect analogy with Horrocks' criterion, I proved that the splitting of a vector 
bundle on any Grassmannian can be verified by restricting it to an arbitrary, standardly 
embedded ${\rm Grass}(2;4)$; this is \emph{far from} being an \emph{ample} 
subvariety. These observations led to consider `sufficiently positive' 
subvarieties of $X$, which posses `many' embedded deformations. 
Loosely speaking, the main result of this article is the following:  
\begin{thm-nono}
Let $\crl V$ be an arbitrary vector bundle on a projective variety $X$, 
defined over an uncountable, algebraically closed field of characteristic zero. 
Let $Y$ be a subvariety of $X$ satisfying the following properties: 
\begin{enumerate}
\item[--] 
$Y$ is $q$-ample for appropriate $q$ 
(e.g. $q+1{\les}\dim X-3\codim_X(Y)$, but not necessarily); 
\item[--] 
The embedded deformations of $Y$ cover an open subset of $X$, 
and their intersection pattern is sufficiently non-trivial. 
\end{enumerate}
\nit Then $\crl V$ splits if and only if it splits along 
a \emph{very general} deformation of $Y$. 
\end{thm-nono}
The precise statement can be found in theorem \ref{thm:split}; applications to the case 
of globally generated vector bundles and of multiplicative group actions, 
with emphasis on homogeneous varieties, are stated in the theorems 
\ref{thm:split-vb1} and \ref{thm:split-homog} respectively. 

A central notion in this article is that of a \emph{$p^\pos/q$-ample subvariety} of 
a projective variety (cf. Definition \ref{def:p>0}); the definition is inspired from 
\cite{ot,to,ar}, also \cite{dps}. The reference \cite{ot} studies in depth the geometric 
properties of the $0$-ample subvarieties. Although the generalization from $0$- to 
$q$-ample is straightforward, it seems that here is the \emph{first place} where 
this weaker property is considered and indeed used for concrete applications. 
A number of geometric properties of $q$-ample subvarieties and also criteria for 
$q$-ampleness, which are necessary for the splitting criterion above, are proved 
in the first section of the article.

The $q$-ampleness condition is asymptotic in nature, involving thickenings 
of subvarieties. For this reason, the price to pay for the generality of the theorem is that 
of testing the splitting along subvarieties which are \emph{very general} within their 
own deformation space. Thus the $q$-ample case studied here can be characterized 
as probabilistic, in contrast with the case of ample vector bundles studied in \cite{hal}, 
which is deterministic. Although the result is formulated algebraically, the proof uses 
complex analytic techniques; the final statement is deduced by base change.

The \emph{essential feature} of the theorem is that of being \emph{intrinsic} 
to the subvariety; it avoids additional data. 
This allows to treat in a \emph{unified way} examples arising 
in totally different situations: 

-- zero loci of globally generated (not necessarily ample) vector bundles 
(cf. sections \ref{sct:0vb}, \ref{sct:split-0loci}), on one hand;

-- homogeneous subvarieties of homogeneous varieties 
(cf. sections \ref{sct:fixed}, \ref{sct:split-homog}), on the other hand. 

\nit As a by-product, we obtain several examples of $q$-ample subvarieties which 
are not zero loci of regular sections in globally generated $q$-ample vector bundles. 

The case of isotropic Grassmannians is detailed in section \ref{sct:grass}. 
As I already mentioned, a vector bundle on $\Grs(u;w)$, 
$u\,{\ges}\,2,\,w\,{\ges}\,u+2$, splits if and only if its restriction to \emph{any} 
embedded $\Grs(2;4)$ splits. Unfortunately, the same proof breaks down in the 
case of the symplectic and orthogonal Grassmannians. 
However, our `probabilistic' approach works very well. 
\begin{thm-nono}
Let $\crl V$ be a vector bundle either on the symplectic Grassmannian $\spGrs(u;w)$, 
with $u\ges 2,w\ges2u$, 
or  the orthogonal Grassmannian $\oGrs(u;w)$, with $u\ges 3,w\ges2u$. 
Then $\crl V$ splits if and only if it splits along a \unbar{\emph{very generally}} 
embedded: \\[1ex]
$\begin{array}[t]{ll}
	\begin{array}{l}-\;\end{array}
	&
	\spGrs(2;4), \text{ in the symplectic case;}
	\\[1ex]  
	\begin{array}{l}-\;\\ \phantom{I}\\ \phantom{I}\end{array}
	&
		\kern-3ex\left.
			\begin{array}{ll}
				\bullet\kern1ex\oGrs(3;6),&\text{for }w=2u,\\ 
				\bullet\kern1ex\oGrs(3;7),&\text{for }w=2u+1,\\ 
				\bullet\kern1ex\oGrs(3;8),&\text{for }w\ges2u+2,
			\end{array}
		\right\}
		\text{ in the orthogonal case.}
\end{array}$ 
\end{thm-nono}
This dramatically simplifies the problem of deciding the splitting of vector bundles 
on Grassmannians. One should compare it with the cohomological criteria 
\cite{ott,ma,ar-ma,ma-oe}, which involve a large number of vanishings.  

Throughout this article, $X$ stands for a smooth, projective, irreducible variety 
over $\bk$, except a few statements involving base change arguments. 
We consider a vector bundle $\crl V$ on $X$ of rank $r$, and denote by 
$\crl E:=\cEnd(\crl V)$ the endomorphisms of $\crl V$. 
For a subscheme $S\subset X$, we let $\crl V_S:=\crl V\otimes\eO_S$, \textit{etc}. 

\setcounter{tocdepth}{1}\tableofcontents


\setcounter{section}{-1}

\section{Prologue: why $q$-ampleness?}\label{sct:why}

Before introducing the relevant definitions, I briefly explain the reasons which led 
to  consider the notion of $p^\pos/q$-ample subvarieties in the context of 
split vector bundles. One says that $\crl V$ splits on $X$ if 
$$
\crl V =\mbox{$\ouset{j\in J}{}{\bigoplus}$}\,\eL_j\otimes\bk^{m_j}
\;\text{with}\;\eL_j\in\Pic(X)\;\text{pairwise non-isomorphic, and}\;
\mbox{$\ouset{j\in J}{}{\sum}$}m_j=r=\rk(\crl V).
$$ 
It is easy to see that $\crl V$ splits if and only if it admits an endomorphism with $r$ 
distinct eigenvalues. The main tool to attack the splitting problem is the following 
observation. 

\begin{m-lemma}\label{lm:idea}
Let $S\subset X$ be a connected, projective subscheme, such that $\crl V_S$ splits. 
If the restriction homomorphism $\res_S:\Gamma(X,\crl E)\to\Gamma(S,\crl E_S)$ 
is surjective, then $\crl V$ splits too. 
\end{m-lemma}
\nit The subschemes considered in this article will usually be thickenings 
of subvarieties of $X$.

\begin{m-definition}\label{def:thick}
For $m\ges 0$, the $m$-th order thickening $Y_m$ of a subvariety $Y\subset X$ is the 
closed subscheme defined by the sheaf of ideals $\eI_Y^{m+1}$; with this convention, 
$Y_0=Y$. The structure sheaves of two consecutive thickenings of $Y$ fit into the exact 
sequence 
\begin{equation}\label{eq:Ym}
0\to\eI_Y^m/\eI_Y^{m+1}\to\eO_{Y_{m}}\to\eO_{Y_{m-1}}\to 0.
\end{equation}
The \emph{formal completion} of $X$ along $Y$ is defined as $\varinjlim Y_m$, 
and is denoted by $\hat X_Y$. 
\end{m-definition}
To apply lemma \ref{lm:idea}, I adopt a two-step strategy:  
\begin{enumerate}\label{strategy}
\item Prove that $\Gamma(X,\crl E)\,{\to}\,\Gamma(\hat X_Y,\crl E_{\hat X_Y})$ 
is surjective; this step requires the $(\dim Y-1)$-ampleness of $Y$. 
(Actually, we will consider also more general thickenings.)

\item Prove that  $\Gamma(\hat X_Y,\crl E_{\hat X_Y})\,{\to}\,\Gamma(Y,\crl E_Y)$ 
is surjective. This step is more involved, requires less ampleness (more positivity) of $Y$, 
and also its genericity in the space of embedded deformations.
\end{enumerate}
If $Y$ has the property that for all vector bundles $\eF$ on $X$ holds 
$H^1(X,\eF\otimes\eI_Y^m)=0$  for all $m$ sufficiently large (depending on $\eF$), 
then the splitting of $\crl V$ along a high order thickening of $Y$ implies 
the splitting on $X$. Indeed, simply take $\eF=\crl E$. 
(For ample subvarieties, explicit lower bounds for $m$ are obtained in \cite{hal}.) 

\begin{m-lemma}\label{lm:y-1}
Assume that $Y$ is $(\dim Y-1)$-ample (cf. \ref{prop:q1} below). Then $\crl V$ splits 
on $X$ if and only if it does so along the formal completion of $X$ along $Y$. 
\end{m-lemma}

The difficulty consists in proving that the splitting of $\crl V$ along a very general 
subvariety $Y\subset X$ implies the splitting of $\crl V$ along $\hat X_Y$.  
This step constitutes the body of the article.


\section{$p^\pos/q$-ample subvarieties}\label{sct:q-tech}

\subsection{Definition and first properties}\label{ssct:q-ample} 

The concept of $q$-ampleness for globally generated vector bundles was introduced in 
\cite{so}; the $q$-ampleness in defined intrinsically in \cite{ar,to} through cohomology 
vanishing properties. We recall the latter definition, the case of globally generated 
vector bundles being detailed in section \ref{sct:0vb}. 

\begin{m-definition}\label{def:q-line}
(cf. \cite[Theorem 7.1]{to}, \cite[Lemma 2.1]{ar}) 
A \emph{line bundle} $\eL$ on a (reduced) projective Gorenstein variety $\tld X$ is 
called \emph{$\tld q$-ample} if for any coherent sheaf $\tld\eF$ on $\tld X$ holds: 
$$
H^t(\tld\eF\otimes\eL^{m})=0,\;\forall\,t>\tld q ,\;\forall\,m\gg0.
$$
We say that a \emph{vector bundle} $\eN$ on $X$ is \emph{$q$-ample} if 
$\eO_{\mbb P(\eN^\vee)}(1)$ on $\mbb P(\eN^\vee)$ is $q$-ample. 
\end{m-definition}

\begin{m-proposition}\label{prop:q1}
{\rm (cf. \cite[Section 7]{to},  \cite[Definition 3.1]{ot}).} 
Let $Y\subset X$ be an equidimensional subscheme of codimension $\delta$, 
$\tld X:=\Bl_Y(X)$ be the blow-up of the ideal of $\,Y$, and $E_Y\subset\tld X$ 
be the exceptional divisor. Assume that $\tld X$ is Gorenstein. 
The following statements are equivalent: 
\begin{enumerate}
\item For any locally free (hence also for any coherent) sheaf $\tld\eF$ on $\tld X$ holds 
\begin{equation}\label{eq:q11}
\begin{array}{r}
H^{t}(\tld X,\tld\eF\otimes\eO_{\tld X}(mE_Y))=0,\;\forall\,t\ges q+\delta,
\;\forall\,m\gg0,
\\ 
(\text{that is, $E_Y$ is $(q+\delta-1)$-ample});
\end{array}
\end{equation}
\item For all locally free sheaves (that is vector bundles) $\eF$ on $X$ holds 
\begin{equation}\label{eq:q12}
H^t(X,\eF\otimes\eI_Y^m)=0,\;\forall\,t\les \dim Y- q,\;\forall\,m\gg0.
\end{equation}
\end{enumerate}
\end{m-proposition}
\begin{proof}
Apply the Serre duality on $\tld X$. 
\end{proof}

If $X$ is smooth and $Y\subset X$ is a locally complete intersection 
(\emph{lci} for short), then $\tld X$ is automatically Gorenstein. 
We are primarily interested in the cohomology vanishing property \eqref{eq:q12}; 
for this reason, we introduce an \textit{ad hoc} terminology. 

\begin{m-definition}\label{def:p>0}
We say that a lci subvariety $Y\subset X$ 
\emph{is (has the property) $p^\pos$} if for any vector bundle $\eF$ on $X$ holds: 
\begin{equation}\label{eq:p}
\exists\,m_{\eF}\ges 1,\kern1ex
\forall\,m\ges m_{\eF},\kern1ex
\forall\,t\les p,\quad
H^t(X,\eF\otimes\eI_Y^{m+1})=0.
\end{equation}
\end{m-definition}
Intuitively, $Y$ is $p^\pos$ if its normal bundle at each point contains at least $p$ 
`positive' directions.
Probably the appropriate name for this property of $Y$ would be 
`\emph{$q$-ample subvariety}', with $q=\dim Y-p$. 
(The case of ample subvarieties \cite{ot} corresponds to $q=0$.) 
The terminology in the definition is made to emphasize the amount of positivity 
of the various objects which appear subsequently. 

\begin{m-proposition}\label{prop:trans}
\nit{\rm(i)} 
 For $Y\subset X$ irreducible, lci holds: 
\begin{equation}\label{eq:equiv-p}
\text{$Y$ is $p^\pos$}
\quad\Leftrightarrow\quad 
\left\{\begin{array}{l}
\text{the normal bundle $\eN=\eN_{Y/X}$ is $(\dim Y-p)$-ample,} 
\\ 
\text{the cohomological dimension $\cd(X\sm Y)\les\dim X-(p+1)$}. 
\end{array}\right.
\end{equation}
(Recall that $\cd(X\sm Y)\ges\codim(Y)-1$ for any closed subvariety.)

\nit{\rm(ii)} 
Let $Z\subset Y$ is $p^\pos$, $Y\subset X$ is $r^\pos$, 
and are both irreducible lci, then 
\begin{equation}\label{eq:trans-p}
\left\{
\begin{array}{l}
\text{$\eN_{Z/X}$ is $\big(\dim Y+\dim Z-(r+p)\big)$-ample,}
\\ 
\text{$\cd(X\sm Z)\les\dim X-(\min\{r,p\}+1)$}.
\end{array}\right.
\end{equation}
In particular,  $Z\subset X$ is $\big(p-(\dim Y-r)\big)^\pos$. 
\end{m-proposition}
For the definition and the properties of the cohomological dimension of a variety, the reader may consult \cite[Ch. III, \S3]{hart-as}.

\begin{proof}
\nit(i) $(\Rightarrow)\;$ Let $\eF$ be a vector bundle on $X$. 
Since $E_Y=\mbb P(\eN)$, the exact sequence 
\begin{equation}\label{eq:mm-1}
0\to\eO_{\tld X}((m-1)E_Y)\to\eO_{\tld X}(mE_Y)\to\eO_{\mbb P(\eN)}(-m)\to 0
\end{equation}
implies $H^t(\mbb P(\eN),\eF\otimes\eO_{\mbb P(\eN)}(-m))=0$, 
for all $t\ges\dim X-p$ and $m\gg0$. 
Let us denote $\delta:=\codim_X(Y)$; the equality 
\begin{equation}\label{eq:PN}
H^{\delta+t}
(\mbb P(\eN),\eO_{\mbb P(\eN)}(-\delta-m)\otimes(\eF\otimes\det(\eN)^{-1}))
= H^{t+1}(\mbb P(\eN^\vee),\eF\otimes\eO_{\mbb P(\eN^\vee)}(m)),
\end{equation}
yields $H^{t+1}(\mbb P(\eN^\vee),\eF\otimes\eO_{\mbb P(\eN^\vee)}(m))=0$, 
for all $t\ges\dim Y-p$. The second implication is \cite[Proposition 5.1]{ot}, 
as $E_Y$ is $\tld q=\delta+(\dim Y-p)-1$ ample. 

\nit$(\Leftarrow)\;$ The $(\dim Y-p)$-ampleness of $\eO_{\mbb P(\eN^\vee)}(1)$, 
combined with \eqref{eq:mm-1}, implies 
\\[.5ex] \centerline{
$H^t(\eF\otimes\eO_{\tld X}((m-1)E_Y))\srel{\cong}{\to}
H^t(\eF\otimes\eO_{\tld X}(mE_Y))$, for $t\ges\dim X-p$ and $m\gg0$,
}\\[1ex] 
so 
$\disp H^t(\eF\otimes\eO_{\tld X}(mE_Y))
=\varinjlim_{k}H^t(\eF\otimes\eO_{\tld X}(kE_Y))$; 
according to \cite[(5.1)]{ot}, the right hand side equals $H^t(X\sm Y,\eF)$, 
which vanishes. 
\smallskip 

\nit(ii) In the sequence 
$$
0\to\eN_{Z/Y}\to\eN_{Z/X}\to{\eN_{Y/X}\big|}_{Z}\to 0,
$$ 
the extremities are $(\dim Z-p)$, respectively $(\dim Y-r)$-ample; 
the subadditivity of the amplitude \cite[Theorem 3.1]{ar} yields the conclusion. 
The bound on the cohomological dimension follows by repeating \emph{ad litteram} 
\cite[Proposition 6.4]{ot}; for the comfort of the reader, we recall the proof here. 
Let $U_Z:=X\sm Z, U_Y:=X\sm Y$ and consider an arbitrary sheaf $\eF$ on $X$. 
In the exact sequence 
$$
\ldots\to H^i_{Y\sm Z}(U_Z,\eF)\to H^i(U_Z,\eF)\to H^i(U_Y,\eF)\to\ldots,
$$
the right hand side vanishes for $i\ges\dim X-r$, because $Y\subset X$ is $r^\pos$. 
We claim that the left hand side vanishes too, for $i\ges \dim X-p$. 
Indeed, it can be computed by using the spectral sequence 
(cf. \cite[Expos\'e I, Th\'eor\`eme 2.6]{groth}): 
$$
H^b(U_Z,\cal H^a_{Y\sm Z}(\eF))\;\Rightarrow\;H^{a+b}_{Y\sm Z}(U_Z,\eF),
$$
where $\cal H^a_{Y\sm Z}(\eF)$ stands for the local cohomology sheaf, with 
support on $Y\sm Z$. The term on the left has the following two properties: 
first, $\cal H^a_{Y\sm Z}(\eF)=\underset{m}{\varinjlim}\,
{\cE}{\kern-1pt}xt^a(\eO_{U_Z}/\eI_{Y\sm Z}^m,\eF)$ 
(cf. \cite[Expos\'e II, Th\'eor\`eme 2]{groth}), the ${\cE}{\kern-1pt}xt$ groups 
are supported on $Y\sm Z$, and $Z\subset Y$ is $p^\pos$, all together imply that 
$$
H^b(U_Z,\cal H^a_{Y\sm Z}(\eF))=0,\;\forall\,b\ges\dim Y-p;
$$
second, $\cal H^a_{Y\sm Z}(\eF)=0,\;\forall\,a\ges\dim X-\dim Y+1,$ because 
$Y\subset X$ is lci. All together, we deduce that 
$H^i_{Y\sm Z}(U_Z,\eF)=0$, for $i\ges(\dim X-\dim Y)+(\dim Y-p-1)+1$. 
For the final statement, observe that $Z$ satisfies the conditions (i). 
\end{proof}


\begin{m-remark}\label{rmk:conn} 
We recall from \cite[Ch. III, Corollary 3.9]{hart-as} the following: 
if $X$ is Gorenstein and $Y\subset X$ is a closed subvariety such that 
$\cd(X\sm Y)\les\dim X-2$, then $Y$ is connected.\footnote{In \textit{loc. cit.} 
the result is stated for $X$ smooth. However, this assumption is used only 
to apply the Serre duality; so the same proof works in the Gorenstein case too. 
We are interested in this generalization for $X$ an lci variety, to avoid 
additional transversality assumptions.} 

In particular, if $Y\subset X$ is lci and $1^\pos$, then it is connected, 
equidimensional. (This is straightforward: take $\eF=\eO_X$ in \eqref{eq:p}.) 
In particular, if $Y$ is smooth then it is irreducible. 
\end{m-remark}

\begin{m-example}\label{expl:strange}
(an unexpected behaviour)\quad 
The example shows that in \eqref{eq:equiv-p} both conditions are necessary. 
For $X:=\mbb P^n\times\mbb P^n$ and $Y:=\text{the diagonal}\cong\mbb P^n$, the following hold:
\begin{enumerate}
\item 
$\euf N_{Y/X}\cong T_{\mbb P^n}$, so the normal bundle is ample 
and globally generated;
\item 
The morphism  
$X\setminus Y\to{\rm Flag}(1, 2;\mbb C^{n+1}),\; 
(x_1,x_2)\mapsto \big(x_1,\lran{x_1,x_2}\big)$, 
is an $\mbb A^1$-fibre bundle, hence $\cd(X\setminus Y)=2n-1$ 
and $Y$ is only $0^{>0}$ in $X$. 
\end{enumerate} 
\end{m-example}

The lack of sufficient positivity of the normal bundle $\eN_{Z/X}$ prevents to conclude 
that $Z\subset X$ is $\min\{r,p\}^\pos$. However, the next proposition shows that it is 
close to be so.

\begin{m-definition}\label{def:ap>0} 
\nit(i) 
For any vector bundle $\eF$ on $X$ and a closed subscheme $Z\subset X$, define 
$$
\tld H^t(Z,\eF) :=
\bigg\{
h\in H^t(Z,\eF)\,\biggl|
\begin{array}{l}
\exists\,\cU\supset Z\text{ open subset of }X,\\ 
\exists\,\tld\alpha\in H^t(\cU,\eF)\text{ such that }\alpha=\res^\cU_{Z}(\tld\alpha)
\end{array}\biggr.
\bigg\}.
$$
Here $\cU$ can be \emph{either} a \emph{Zariski or} an \emph{analytic} 
(tubular) open neighbourhood of $Z$. 

\nit(ii) 
We say that a subscheme $Z\subset X$ is $p^\apos$ (\emph{approximately} $p^\pos$) 
if there is a decreasing sequence of sheaves of ideals $\{\eJ_m\}_{m}$ such that 
the following hold:  
\begin{equation}\label{eq:ap}
\begin{array}{ll}
\bullet\;&
\forall\,m,n\ges1\;\exists\,m'>m,\,n'>n\text{ such that }
\eJ_{m'}\subset\eI_Y^m,\;\eI_Y^{n'}\subset\eJ_{n};
\\[1ex] 
\bullet\;&
\text{for any vector bundle $\eF$ on $X$,}
\\[.5ex]
&
\exists\,m_\eF\ges1\text{ such that }H^t\big(X,\eF\otimes{\eJ_{m}}\big)=0,
\;\forall\,m\ges m_\eF\;\forall\,t\les p.
\end{array}
\end{equation}
\end{m-definition}

\begin{m-proposition}\label{prop:approx-p}
If $Z\subset Y,\;Y\subset X$ be lci, and both $p^\pos$, then $Z\subset X$ is $p^\apos$. 
In particular, for all vector bundles $\eF$ on $X$ holds: 
\\[1ex] \centerline{
$
\res^X_{Z}:H^t(X,\eF)\to\tld H^t(Z,\eF)\text{ is: }
\Big\{\begin{array}{l}
\text{- an isomorphism, for }t\les p-1,\\ 
\text{- injective, for }t=p.
\end{array}\Big.
$
}
\end{m-proposition}

\begin{proof}
Since $Z\subset Y\subset X$ are lci, for any $l\ges a$ one has the exact sequence:
\begin{equation}\label{eq:la}
0\to
\frac{\eI_Y^a}{\eI_Y^{a+1}}\otimes\biggl(\frac{\eI_Z}{\eI_Y}\biggr)^{l-a}
\to
\frac{\eI_Z^l+\eI_Y^{a+1}}{\eI_Y^{a+1}}
\to
\frac{\eI_Z^l+\eI_Y^{a}}{\eI_Y^{a}}
\to0.
\end{equation}
Indeed, just use 
$(\eI_Z^l+\eI_Y^{a})/{\eI_Y^{a}}\cong\eI_Z^{l}/\eI_Y^a\cdot\eI_Z^{l-a}$. 
Note that the left hand side is an $\eO_Y$-module; 
$\eI_Z/\eI_Y=\eI_{Z\subset Y}$ is the ideal of $Z\subset Y$, and 
that $\eI_Y^a/\eI_Y^{a+1}=\Sym^a\eN_{Y/X}^\vee$ is the symmetric power 
of co-normal bundle. 
The $p^\pos$ assumption implies that there are integers $k_\eF,l_\eF$, 
and a linear function $l(k)=\l k+\mu$ ($\l,\mu$ independent of $\eF$) 
with the following properties: 
\\[1ex] \centerline{
$
\begin{array}{rl}
H^t(\eF\otimes\eI_{Y}^k)=0,
&
\;\forall\,t\les p,\;\forall\,k\ges k_\eF,
\\[.5ex]
H^t(\eF_Y\otimes\eI_{Z\subset Y}^{l})=0,
&
\;\forall\,t\les p,\;\forall\,l\ges l_\eF,
\\[.5ex] 
H^t(\eF_Y\otimes\Sym^a\eN_{Y/X}^\vee\otimes\eI_{Z\subset Y}^{l-a})=0,
&
\;\forall\,t\les p,\;\forall\,a\les k,\;\forall\,l\ges l(k).
\end{array}
$
}\\[1ex]
The last claim uses the uniform $q$-ampleness property \cite[Theorem 6.4]{to} 
and the subadditivity of the regularity \cite[Theorem 3.4]{to}: first, 
there is a linear function $l(r)$ such that for any vector bundle $\eF$ 
with regularity ${\rm reg}(\eF_Y)\les r$ holds 
$$
H^t(\eF_Y\otimes\eI_{Z\subset Y}^l)=0,\;\text{for}\;t\les p,\; l\ges l(r);
$$
second, for $a\les k$, 
${\rm reg}(\eF_Y\otimes\Sym^a\eN_{Y/X}^\vee)\les{\rm linear}(k)$. 
Recursively for $a=1,\ldots,k$, starting with 
$\frac{\eI_Z^l+\eI_Y}{\eI_Y}=\eI_{Z\subset Y}^l$, 
the exact sequence \eqref{eq:la} yields: 
\\[.5ex] \centerline{
$
H^t\Big(
\eF\otimes\frac{\eI_Z^l+\eI_Y^{k}}{\eI_Y^{k}}
\Big)=0,\;\forall\,l\ges l(k). 
$
}\\[.5ex]
Now plug this into 
$0\to\eI_Y^k\to\eI_Z^l+\eI_Y^{k}\to\frac{\eI_Z^l+\eI_Y^{k}}{\eI_Y^{k}}\to0$ 
(tensored by $\eF$), and deduce that 
$$
H^t\big(\eF\otimes(\eI_Y^k+\eI_Z^l)\big)=0,
\;\forall\,t\les p,\;\;\forall\,k\ges k_\eF,\;\forall\,l\ges l(k).
$$
We denote by $Z_{l,k}$ the subscheme defined by $\eI_Y^k+\eI_Z^l$;  
it is an `asymmetric' thickening of $Z$ in $X$. 
For any $l$ as above, there is $m>l$ such that  
$$
\eI_Y^l+\eI_Z^m\subset\eI_Z^l\subset\eI_Y^k+\eI_Z^l
\;\Rightarrow\;Z_{l,k}\subset Z_l\subset Z_{m,l}.
$$
For $m>l>k$ as above, one has the commutative diagram 
\\[.5ex] \centerline{
$
\xymatrix@R=.25em@C=4em{
&H^t(\eF_{Z_{m,l}})\ar[dr]^-{\res^{m,l}_l}&
\\ 
H^t(\eF)
\ar[rr]|{\ \xi\ }
\ar[ur]^-{\res^X_{{m,l}}}
\ar[dr]_-{\res^X_{{l,k}}}
&& H^t(\eF_{Z_l}).\ar[dl]^-{\res^l_{l,k}}
\\ 
&H^t(\eF_{Z_{l,k}})&
}
$
}\\[.5ex]
Notice that $\res^l_{l,k}\circ \xi=\res^X_{l,k}$, which is injective, so $\xi$ is injective. 
It remains to prove that $\xi$ maps onto $\tld H^t(\eF_{Z_l})$, for $t\les p-1$: 
assume $\alpha=\res^\cU_{Y_l}(\tld\alpha)$, with $\cU\supset Z$ open and 
$\tld\alpha\in H^t(\eF_\cU)$; then $\alpha=\res^m_l(\,\res^\cU_{{m,l}}(\tld\alpha)\,)$, 
and $\res^\cU_{{m,l}}(\tld\alpha)$ belongs to the image of $\res^{X}_{{m,l}}$.
\end{proof}

\begin{m-lemma}\label{lm:pull-back}
Let $\vphi:X\to X'$ be a flat, surjective morphism, with $X,X'$ smooth, 
whose fibres are $d$-dimensional. If $Y'\subset X'$ is lci and $p^\pos$, 
then $Y=\vphi^{-1}(Y')\subset X$ is the same. 
\end{m-lemma}

\begin{proof}
Since $f$ is flat, $Y$ is still lci in $X$ and $\codim_X(Y)=\codim_{X'}(Y')=\delta$. 
Let us check the property \eqref{eq:q11}. The universality property of the blow-up 
(cf. \cite[Ch. II, Corollary 7.15]{hart-ag})  yields the commutative diagram 
$$
\xymatrix@R=1.5em@C=4em{
\tld X=\Bl_{Y}(X)\ar[d]\ar[r]^-{\tld\vphi}&\tld X'=\Bl_{Y'}(X)\ar[d]
\\
X\ar[r]^-\vphi&X'.
}
$$
The morphism $\tld\vphi$ still has $d$-dimensional fibres and 
$\tld\vphi^*\eO_{\tld X'}(E_{Y'})=\eO_{\tld X}(E_Y)$. 
For any coherent sheaf $\tld\eF$ on $\tld X$ holds: 
$$
R^i\tld\vphi_*(\tld\eF\otimes\eO_{\tld X}(mE_Y))
=R^i\vphi_*\tld\eF\otimes\eO_{\tld X'}(mE_{Y'}), \quad 
R^{>d}\tld\vphi_*\tld\eF=0.
$$ 
As $Y'\subset X'$ is $p^\pos$, we deduce that 
\\[1ex] \centerline{
$H^t\big(\tld X',R^i\tld\vphi_*\tld\eF\otimes\eO_{\tld X'}(mE_{Y'})\big)=0$, 
for $i=0,\ldots,d$, $t\ges(\dim X'-p)$, and $m\gg0$. 
}\\[1ex]
Then the Leray spectral sequence implies that $E_Y$ is 
$\big((\dim X'-p)+d-1\big)$-ample. 
\end{proof}


\subsection{Criterion for the positivity of a subvariety}\label{ssct:+}

Many examples of $q$-ample subvarieties occur as zero loci of regular sections in 
(Sommese) $q$-ample vector bundles. 
(This will be detailed in the section \ref{sct:0vb}.)  However, it will be necessary to test 
the $q$-ampleness of a subvariety in more general circumstances. 

\begin{m-proposition}\label{prop:x-b}
Let the situation be as in \ref{prop:q1}. Assume that there is an irreducible variety $V$ 
and a morphism $b:\tld X\to V$ such that $\eO_{\tld X}(E_Y)$ is $b$-relatively ample. Then $Y\subset X$ is $p^\pos$, for $p:=\dim X-\dim b(\tld X)-1$.\\ 
If $\tld X$ is the blow-up of an ideal $\eJ\subset\eO_X$ with $\sqrt{\eJ}=\eI_Y$ 
and $b$ is as above, then $Y$ is $p^\apos$. 
\end{m-proposition}

\begin{proof}
Let $\tld\eF$ be a coherent sheaf on $\tld X$, and $j\ges\dim X-p>\dim b(\tld X)$. Then 
\\[.5ex] \centerline{
$R^tb_*(\tld\eF\otimes\eO_{\tld X}(mE_Y))=0$, for $t>0$ and $m\gg0$,
}\\[.5ex] 
implies that 
$H^{j}\big(\tld X,\tld\eF\otimes\eO_{\tld X}(mE_Y)\big)=
H^j\big(\,V,b_*(\tld\eF\otimes\eO_{\tld X}(mE_Y))\,\big).$ 
But the right hand side vanishes, because $\Supp b_*(\tld\eF\otimes\eO_{\tld X}(mE_Y))$ 
is at most $\dim b(\tld X)$-dimensional. 
\end{proof}

\begin{m-remark}\label{rmk:x-b}
In the section \ref{sct:split-crit}, we will be interested in families of $p^\pos$ 
subvarieties of $X$; the proposition generalizes straightforwardly. 
Let $\cY\subset S\times X$ be a smooth family of subvarieties of $X$ indexed by 
some parameter space $S$. Assume that there is an $S$-variety $V$ such that  
\\[1ex]\centerline{
$\eO(E_\cY)$ is relatively ample for a morphism $\Bl_{\cY}(S\times X)\to V$. 
}\\[1ex]
Then $Y_s$ is $(\dim (S\times X)-\dim V-1)^\pos$, for all $s\in S$ such that $Y_s$ is lci.
\end{m-remark}

The criterion applies to the situations discussed in sections \ref{sct:0vb} and 
\ref{sct:fixed}, that is the case of zero sections in globally generated vector bundles, 
and of sources of $G_m$-actions.


\subsection{$q$-positive line bundles}\label{ssct:q-line-bdl}

Intuitively, the hypotheses in \ref{prop:x-b} allow $\eO_{\tld X}(E_Y)$ to be `negative' 
in at most $\dim b(\tld X)<\dim X-p$ directions. 
As we will see below, this stronger positivity property allows to control $\Pic(Y)$. 

\begin{m-definition}\label{def:p-line-bdl}{(cf. \cite{dps})} 
A line bundle $\eL$ on a smooth projective variety $\tld X$ is \emph{$q$-positive} 
if it admits a Hermitian metric whose curvature form has at most $q$ negative (or zero) 
eigenvalues. 
\end{m-definition}
It is generally true (cf. \cite[Proposition 28]{ag}, \cite[Proposition 2.1]{dps}) 
that a $q$-positive line bundle $\eL$ satisfies $H^{\ges q+1}(\tld\eF\otimes\eL^m)=0$, 
for any vector bundle $\tld\eF$ and $m\gg0$, that is 
\\[1ex] \centerline{
$\eL$ is $q$-positive
$\;\Rightarrow\;$ 
$\eL$ is $q$-ample. 
}

\begin{m-proposition}\label{prop:+}
Let the situation be as in proposition \ref{prop:x-b}, with $X,Y$ smooth. 
Then $\eO_{\tld X}(E_Y)$ is $\dim b(\tld X)$-positive. 
\end{m-proposition}

\begin{proof}
One may assume that $V$ is smooth; otherwise, take an embedding into a smooth variety. 
As $\eO_{\tld X}(E_Y)$ is $b$-relatively ample, there is an embedding 
$\tld X\srel{\iota}{\to}\mbb P^N\times V$ (over $V$) such that 
$\eO_{\tld X}(m_0E_Y)=\iota^*(\eO_{\mbb P^N}(1)\boxtimes\euf M)$, 
for some $m_0>0$ and $\euf M\in\Pic(V)$. 
Take a strictly positive Hermitian metric on $\eO_{\mbb P^N}(1)$, an arbitrary 
on $\euf M$, and the product metric on $\eO_{\mbb P^N}(1)\boxtimes\euf M$. 

The maximal rank of $\rd b_{\tld x}$ is $\dim b(\tld X)$, attained on a dense open 
subset of $\tld X$, so 
$$
\rk(\rd b_{\tld x})\les\dim b(\tld X),\;\forall\,\tld x\in\tld X.
$$ 
Since $\iota$ is an embedding, the curvature of the pull-back metric on 
$\eO_{\tld X}(m_0E_Y)$ is positive definite on $\Ker(\rd b_{\tld x})$, and 
$\dim\Ker(\rd b_{\tld x})\ges\dim X-\dim b(\tld X)$. 
So, at each point of $\tld X$, there are at most $\dim b(\tld X)$ negative eigenvectors. 
\end{proof}

\begin{m-lemma}\label{lm:pull-back-smooth}
Let $\vphi:X\to X'$ be a smooth morphism of relative dimension $d$. If $\eL'$ is 
a $q$-positive line bundle on $X'$, then $\vphi^*\eL'$ is $(q+d)$-positive on $X$. 
\end{m-lemma}

\begin{proof}
Take a metric in $\eL'$ as in \ref{def:p-line-bdl}, and pull it back to $\eL$. 
\end{proof}


\subsection{Geometric properties of $q$-ample subvarieties}\label{ssct:geom}

In the section \ref{sct:split-crit} we will be concerned with families of $p^\pos$ subvarieties. 
The non-emptiness and connectedness of their intersections, as well as 
their Picard groups will be essential; these issues are investigated below. 


\subsubsection{Non-emptiness and connectedness of the intersections}\label{sssct:conn} 

For a family of irreducible subvarieties $\cY\subset S\times X$ of $X$ indexed by some 
variety $S$, we denote 
\begin{equation}\label{eq:ost}
Y_{st}:=Y_s\cap Y_t,\quad Y_{ost}:=Y_o\cap Y_s\cap Y_t,\;\forall\,o,s,t\in S.
\end{equation}

\begin{m-lemma}\label{lm:yy'} 
Let $\cY=\{Y_s\}_{s\in S}$ be a family of $\delta$-codimensional, lci subvarieties 
of $X$ such that their double and triple intersections are equidimensional lci 
(if non-empty). Assume that 
$$
\cd(X\sm Y_s)\les\dim(X)-2\delta-2,\;\forall\,s\in S.
$$
Then, for all $o,s,t\in S$, the intersections $Y_{os}$ and $Y_{ost}$ are indeed 
\emph{non-empty} and \emph{connected}.  
In particular, the conclusion holds if $\cY$ is a family 
of $p^\pos$ subvarieties of $X$, with $2\delta+1\les p$.
\end{m-lemma}

\begin{proof}
If $Y_{os}=\emptyset$, then $Y_s\subset X\sm Y_o$ implies 
$\dim X-2\delta-2\ges\cd(X\sm Y_o)\ges\dim X-\delta,$ a contradiction. 
Hence the double intersections are non-empty. 
If $Y_{ost}=\emptyset$, then $Y_{os}\subset Y_o\sm Y_{ot}$, hence 
we obtain $\dim X-2\delta-2\ges\cd(Y_o\sm Y_{ot})\ges\dim X-2\delta,$ 
another contradiction. 
The connectedness of the intersections follows from the remark \ref{rmk:conn}: 
$$
\begin{array}{l}
\cd(Y_o\sm Y_{os})\les\cd(X\sm Y_s)\les\dim X-2\delta-2<\dim Y_o-2;
\\ 
\cd(Y_{os}\sm Y_{ost})\les\cd(X\sm Y_t)\les\dim X-2\delta-2\les\dim Y_{os}-2.
\end{array}
$$
The last claim follows from the proposition \ref{prop:trans}.  
\end{proof}


\subsubsection{The Picard group of $p^\pos$ subvarieties}\label{sssct:pic}

Let $Y\subset X$ be smooth. We are interested when the pull-back 
$\res^X_Y:\Pic(X)\to\Pic(Y)$ is an isomorphism. 

\begin{m-lemma}\label{lm:pic-iso}
Let $Y\subset X$ be a subvariety. 
Then $\res^X_Y:\Pic(X)\to\Pic(Y)$ is an isomorphism as soon as 
$H^t(X,\mbb Z)\to H^t(Y;\mbb Z),\; t=1,2,$ are isomorphisms. 
In particular, $\res^X_Y$ is an isomorphism if 
$H_{2\dim_{\bk}X-t}(X\sm Y;\mbb Z)=0$ for $1\les t\les3$. 
\end{m-lemma}

\begin{proof}
The hypothesis implies $H^t(X;\bk)\srel{\cong}{\to}H^t(X;\bk)$, for $t=1,2$; 
the Hodge decomposition yields $H^t(X;\eO_X)\srel{\cong}{\to}H^t(X;\eO_Y)$. 
Now compare the exponential sequences for $X$ and $Y$: 
\begin{equation}\label{eq:h12pic}
\xymatrix@R=1.25em{
H^1(X;\mbb Z)\ar[r]\ar[d]^-{\res^X_Y}&
H^1(X;\eO_X)\ar[r]\ar[d]&
\Pic(X)\ar[r]\ar[d]&
H^2(X;\mbb Z)\ar[r]\ar[d]^-{\res^X_Y}&
H^2(X;\eO_X)\ar[d]
\\ 
H^1(Y;\mbb Z)\ar[r]&
H^1(Y;\eO_X)\ar[r]&
\Pic(Y)\ar[r]&
H^2(Y;\mbb Z)\ar[r]&
H^2(Y;\eO_X).
}
\end{equation}
\end{proof}

\begin{m-remark}\label{rmk:q-ample}
The restriction $H^t(X;\bbQ)\to H^t(Y;\bbQ)$ is an isomorphism for $t\les p-1$, 
for any $p^\pos$ subvariety $Y$ (cf. \cite[Corollary 5.2]{ot}).  
Hence, if $Y$ is $3^\pos$, $\Pic^0(X)\to\Pic^0(Y)$ is a finite morphism, 
and $\NS(X)\to\NS(Y)$ has finite index. 
\end{m-remark}

\begin{m-theorem}\label{thm:pic}
Let $Y\subset X$ be a smooth subvariety. Assume there is a smooth variety $V$ 
and a morphism $\tld X=\Bl_Y(X)\srel{b}{\to}V$ such that $\eO_{\tld X}(E_Y)$ is 
$b$-relatively ample. Then holds 
\\[.5ex] \centerline{
$H^t(X;\mbb Z)\srel{\cong\,}{\to}H^t(Y;\mbb Z)$, 
for $t\les p-1$ and  $p:=\dim X-\dim b(\tld X)-1$.
}\\[1ex] 
Thus $\Pic(X)\to\Pic(Y)$ is an isomorphism, for $p\ges 3$. 
If $Y'$ is a small, smooth deformation of $Y$, then $\Pic(X)\to\Pic(Y')$ 
is still an isomorphism, and $Y'$ is $p^\pos$ too.
\end{m-theorem}

\begin{proof}
The hypothesis implies that $\eO_{\tld X}(E_Y)$ is $\dim b(\tld X)$-positive 
(cf. proposition \ref{prop:+}). Then \cite[Theorem III$\,$]{bott-lefschetz}, 
\cite[Lemma 10.1]{ot} imply that 
$$
H^t_c(X\sm Y;\mbb Z)=H^t_c(\tld X\sm E_Y;\mbb Z)=0,
\;\forall\,t\les\dim X-\dim b(\tld X)-1.
$$
The long exact sequence in cohomology yields the conclusion. 
The second statement is the openness of the $q$-ample property \cite[Theorem 8.1]{to}. 
\end{proof}


\section{The splitting criterion}\label{sct:split-crit}

Let $\crl V =\ouset{j\in J}{}{{\mbox{$\bigoplus$}}}\,\eL_j\otimes\bk^{m_j}$ 
be a split vector bundle, with $\eL_j\in\Pic(X)$ pairwise non-isomorphic. 
In this case, $\crl V_j:=\eL_j\otimes\bk^{m_j}$, ${j\in J}$, are called the 
{\em isotypical components} of the splitting. 
If $\ouset{j\in J}{}{{\mbox{$\bigoplus$}}}\,\eL_j\otimes\bk^{m_j}$ and 
$\ouset{j\,'\in J'}{}{{\mbox{$\bigoplus$}}}\,\eL'_{j\,'}\otimes\bk^{m\,'_{j\,'}}$ are 
two splittings, then there is a bijective function $J\srel{\epsilon}{\to} J'$ 
such that $\eL'_{\epsilon(j)}\cong\eL_j$ and $m'_{\epsilon(j)}=m_j$ for all $j\in J$ 
(cf. \cite[Theorem 1 and 2]{at}.) However, the isotypical components are 
{\em not uniquely defined}, because the global automorphisms of $\crl V$ send a 
splitting into a new one. We consider the partial order `$\prec$' on $\Pic(X)$: 
\begin{equation}\label{eq:order}
\eL\prec\eL'\quad\text{if }\eL\neq\eL'\text{ and }\Gamma(X,\eL^{-1}\eL')\neq 0.
\end{equation}
The isotypical components corresponding to the 
\emph{maximal} elements, are \emph{canonically defined}. 

\begin{m-lemma}\label{lm:max}
Let $M\subset J$ be the subset of maximal elements with respect to $\prec$. 
Then there is a natural, injective homomorphism of vector bundles 
\begin{equation}\label{eq:evm}
\ev_M:\underset{j\in M}{\mbox{$\bigoplus$}}
\eL_j\otimes\Gamma(X,\eL_j^{-1}\otimes\crl V)
\to \crl V.
\end{equation}
\end{m-lemma}
\nit This is an essential observation because, 
given a family of subvarieties $\{Y_s\}_{s\in S}$ of $X$ such that $\crl V$ splits along 
each of them, we can \textit{glue together} the maximal isotypical components of 
$\crl V_{Y_s}$ (cf. lemma \ref{lm:glue} below). 

\subsection{Gluing of split vector bundles}\label{ssct:glue}

Let $S$ be an irreducible quasi-affine variety,  and $\cY\subset S\times X$ be subvariety 
(denote by $\pi,\rho$ the morphisms to $S$ and $X$, respectively) 
with the following properties: 

\begin{longtable}%
{>{\raggedleft}p{3.5ex}%
>{\raggedright}p{0.775\linewidth}|%
>{\raggedleft}p{0.085\linewidth}}
	(0)
&
	$\Pic(S)$ is trivial. (This can be achieved by shrinking $S$.) 
&
\tabularnewline 
	(i)
&
	The morphism $\cY\srel{\pi}{\to} S$ is proper, smooth, with connected fibres. 
&
\tabularnewline 
	(ii)
&
	For $o\in S$, define $S(o):=\{s\in S\mid Y_s\cap Y_o\neq\emptyset\}$. 
	(It is a subscheme of $S$.) \\ 
	Assume that the \emph{very general} point $o\in S$  has the following properties: 
&
\tabularnewline 
&
	(a) $Y_{os}=Y_o\cap Y_s$ is irreducible, for all $s\in S(o)$,  
	and the intersection is transverse for $s\in S(o)$ generic;
&
\tabularnewline 
&
	(b) 
	For $s\in S(o)$ generic, all the arrows below are isomorphisms: 
&
\tabularnewline 
&
	\xymatrix@R=.75em@C=0.12\linewidth{
	&
	&
		\Pic(Y_s)\ar[dr]^-{\res^{Y_s}_{Y_{os}}}
	&
	&
	\\ 
	&
		\Pic(X)\ar[ur]^-{\res^X_{Y_s}}
		\ar[dr]_-{\res^X_{Y_o}}\ar[rr]|{\;\res^X_{Y_{os}}\;}
	&
	&
		\Pic(Y_{os});
	&
		(\Pic)
	\\ 
	&
	&
		\Pic(Y_o)\ar[ur]_-{\res^{Y_o}_{Y_{os}}}
	&
	&
	}
&
\tabularnewline 
&
	(c) 
	$\rho(\cY_{S(o)})\subset X$ is open:
&
	\begin{minipage}[b]{3.5em}
		\ysx\hskip-1ex
	\end{minipage}
\tabularnewline 
&
		$\begin{array}{cclcc}
			\phantom{o}
		&
			\xymatrix@R=1.5em@C=0.05\linewidth{
				\ar@{-}[rrddd]_-{Y_s}&&&
			\\ 
				\ar@{-}[rrr]^-{Y_o}&&&
			\\ 
				&\ar@{}[r]_-\bullet&&
			\\ 
				&\ar@{}[r]^{{\raise1ex\hbox{$x\phantom{xi}$}}}&&
			}
		&
			\quad
		&
			\begin{minipage}[t]{30ex}\vskip2ex 
			The \arm condition:\\ these configurations cover\\ an open neighbourhood of $Y_o$.
			\end{minipage}
		&
			\kern0.07\linewidth\lower6ex\hbox{(1\text{-arm})}
		\end{array}$
& 
\tabularnewline 
&
	(d) 
	For all $s,t\in S(o)$ holds: 
&
\\ 
&
	$\begin{array}{clr}
		\phantom{o}
	&
		Y_{st}\neq\emptyset\;\Rightarrow\;Y_{ost}\neq\emptyset\quad\text{and connected}.
	&
	\\ 
	&
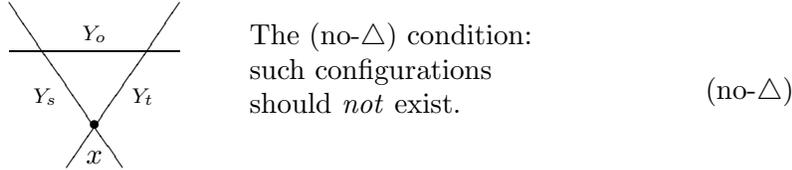
		\begin{array}{ccl}
			\xymatrix@R=1.5em@C=0.05\linewidth{
				\ar@{-}[rrddd]_-{Y_s}
			&
			&
			&\ar@{-}[llddd]^-{Y_t}
			\\ 
				\ar@{-}[rrr]^-{Y_o}
			&
			&
			&
			\\ 
			&\ar@{}[r]_-\bullet
			&
			&
			\\ 
			&
			\ar@{}[r]^{{\hbox{$x$}}}
			&
			&
			}
		&
			\quad
		&
			\begin{minipage}[t]{25ex}\vskip2ex 
			The \notr condition:\\ such configurations\\ should \emph{not} exist.
			\end{minipage}
		\end{array}
	&
	\kern0.115\linewidth\lower1ex\hbox{(\text{no-}$\triangle$)}
	\end{array}$
&
\end{longtable}

\begin{m-remark}\label{rmk:pic}
If all the triple intersections are non-empty and connected 
(sufficient conditions are given in \ref{lm:yy'}), 
then $S(o)=S$, and \arm\kern-2pt, \notr are automatically satisfied. 
This property holds in the case of zero loci of sections in ample vector bundles, 
which was studied in \cite{hal}. More general situations when \ysx 
is satisfied appear in \ref{prop:pic-0loci}, \ref{thm:pic-sources}. 

However, for the isotropic symplectic and orthogonal Grassmannian studied 
in section \ref{sct:grass}, the generic double and triple intersections are empty. 
The analysis of these cases led to the weaker conditions above. 
\end{m-remark}

\begin{m-definition}\label{def:ggen}
(i) The \emph{geometric generic} fibre of $\pi$ is defined by the Cartesian diagram 
\\ \centerline{
$
\xymatrix@R=1.25em@C=2em{
\mbb Y\ar[r]\ar[d]&\cY\ar[d]^-\pi&
\\ \Spec(\bar\Bbbk)\ar[r]&S,& \Bbbk:=\bk(S).
}
$
}
A \emph{very general point} of $S$ refers to a point outside a countable union 
of subvarieties of $S$. 

\nit(ii) 
The \emph{double} (respectively \emph{triple}) \emph{geometric generic self-intersections} 
of $\cY$, denoted $\mbb Y_2$ and $\mbb Y_3$ respectively, are defined by means of 
the diagrams: 
\begin{equation}\label{eq:y23}
\begin{array}{lc}
\xymatrix@R=1.25em@C=1.2em{
\cY_2:=\overline{(\cY\times_{X}\cY)\sm\diag(\cY)}\ar@{^(->}[r]\ar[dd]^-{\pi_2}&
\cY\times_X\cY\ar[r]\ar@{^(->}[d]
&
X\ar@{^(->}[d]^{\text{diag}}
\\ 
&
\cY\times\cY\ar[r]\ar[d]
&
X\times X
\\ 
S_2:=\Supp(\pi\times\pi)_*\eO_{\cY_2}
\ar@{^(->}[r]&S\times S&
}
&\quad
\xymatrix@R=1.5em@C=1.2em{
\cY_3\hra\cY\times_{X}\cY\times_X\cY\ar[d]\ar[r]&X
\\ 
S\times S\times S.
}
\\ & \\ 
\mbb Y_2:=\cY_2\times\Spec\big(\,\ovl{\bk(S_2)\!}\;\big),&
\mbb Y_3\text{ is defined similarly.}
\end{array}
\end{equation}
Thus there are natural morphisms $\mbb Y_2 \rightrightarrows\mbb Y$ 
(the projections onto the first and second components), and 
$\mbb Y_3\,\raise6.5pt\hbox{$\rightarrow$}\kern-14pt\rightrightarrows\mbb Y_2$. 
(One may think off $\bbY_2$ as $Y_{st}=Y_s\cap Y_t$, for very general $(s,t)\in S_2$, 
and the morphisms are the inclusions into $Y_s, Y_t$ respectively; similarly for 
$\bbY_3$.) 
\end{m-definition}
For $o\in S$, $S(o)$ introduced at \ysx is 
$(\pr^{S_2}_S)^{-1}(o)\subset(\{o\}\times S)\cap S_2$.

\begin{m-lemma}\label{lm:ball}
There is an open analytic subset (a ball) $\ball\subset S$, which can be chosen arbitrarily 
in some Zariski open subset of $S$, such that the following hold: 

\nit{\rm (i)} If $(\rho^*\crl V)\otimes\bar\Bbbk$ splits on $\bbY$, then 
$(\rho^*\crl V)_\ball$ splits.  

\nit{\rm (ii)} If $\Pic(X\otimes_{\bk}\bar\Bbbk)\to\Pic(\mbb Y)$ and 
$\Pic(\mbb Y) \rightrightarrows \Pic(\mbb Y_2)$ are isomorphisms, 
then \ysx\hskip-1ex$(\Pic)$ is satisfied for points in $\ball$. \newline 
This condition is satisfied in the particular case when $\cY\subset S\times X$ 
is a $\delta$-codimensional, $p^\pos$ family of subvarieties of $X$ 
as in \ref{rmk:x-b}, with $p\ges\delta+3$.  
\end{m-lemma}
 
\begin{proof}
Since $(\rho^*\crl V)_\bbY$ splits, there are $\ell'_1,\ldots,\ell'_r\in\Pic(\bbY)$ 
such that $(\rho^*\crl V)_\bbY=\ell'_1\oplus\ldots\oplus\ell'_r$; 
the right hand side is defined over a finitely generated (algebraic) extension of $\bk(S)$. 
Thus there is an open affine $S^\circ\subset S$, and a finite morphism $\si:S'\to S^\circ$ 
such that $\ell'_1,\ldots,\ell'_r$ are defined over $\bk[S']$, and $(\rho^*\crl V)_{S'}$ 
splits on $\cY_{S'}$. 
Then there are open balls $\ball'\subset S'$ and $\ball\subset S^\circ$ such that 
$\ball'\srel{\si}{\to}\ball$ is an {\em analytic} isomorphism, and the splitting of 
${(\rho^*\crl V)}_{\ball'}$ on $\cY_{\ball'}$ descends to ${(\rho^*\crl V)}_{\ball}$ 
on $\cY_{\ball}$. The second statement is analogous. 

In the situation \ref{rmk:x-b}, $\Pic(X_{\bar\Bbbk})\to\Pic(\bbY)$ 
is an isomorphism for $p\ges 3$, cf. theorem \ref{thm:pic}; 
the isomorphism $\Pic(\bbY)\to\Pic(\bbY_2)$, requires $p-\delta\ges 3$. 
\end{proof}

Henceforth, for $o\in S$ very general, we replace $S$ by $S(o)$ 
(actually by an irreducible component of it containing $o$) 
and $\cY$ by $\cY\times_S{S(o)}$; 
the restrictions of $\pi,\rho$ are denoted the same. 

\begin{m-lemma}\label{lm:rho}
$\Pic(X)\srel{\rho^*}{\to}\Pic(\cY_{S})$ is an isomorphism, for $S$ small enough. 
\end{m-lemma}

\begin{proof}
The composition $\Pic(X)\srel{q^*}{\to}\Pic(\cY)\srel{\res_{Y_o}}{\to}\Pic(Y_o)$ 
is bijective, so $\rho^*$ is injective. 
For the surjectivity, take $\ell\in\Pic(\cY)$. If $\ell_{Y_o}\cong\eO_{Y_o}$, then 
$$
\{s\in S\mid \ell_{Y_s}\not\cong\eO_{Y_s}\}=\{s\in S\mid h^0(\ell_{Y_s})=0\}
$$
is open, by semi-continuity, so $\{s\in S\mid\ell_{Y_s}\cong\eO_{Y_s}\}$ is closed. 
On the other hand, by restricting to $Y_{os}$, the hypothesis $(\Pic)$ implies that this 
set is dense; thus it is the whole $S$. 
It follows that $\ell\cong\pi^*\bar\ell$, with $\bar\ell\in\Pic(S)$. 
But $\Pic(S)$ is trivial for $S$ sufficiently small, so $\ell\cong\eO$. 
If $\ell\in\Pic(\cY)$ is arbitrary, take $\eL\in\Pic(X)$ such that 
$\ell_{Y_o}\cong\eL_{Y_o}$, so ${(\rho^*\eL^{-1})\ell|}_{Y_o}$ is trivial. 
\end{proof}

Let $\crl V$ be a vector bundle on $X$, such that $\rho^*\crl V$ splits; 
by the previous lemma, 
\begin{equation}\label{eq:VL}
{\rho^*\crl V}\cong\rho^*\bigl(\underset{j\in J}{\oplus}\eL_j^{\oplus m_j}\bigr),
\text{ with $\eL_j\in\Pic(X)$ pairwise non-isomorphic.}
\end{equation}
For any $s\in S$, let $M_s\subset J$ be the subset of maximal elements, 
for $\crl V_{Y_s}$. 
By semi-continuity, there is a neighbourhood $S_s\subset S$ of $s$ such that 
$M_s\subset M_{s'}$ for all $s'\in S_s$; thus there is a largest subset $M\subset J$, 
and an open subset $S'\subset S$ such that $M=M_s$, for all $s\in S'$. 
Hence, possibly after shrinking $S$, the set of maximal isotypical components of 
$\crl V_{Y_s}$ with respect to \eqref{eq:order} is independent of $s\in S$.

\begin{m-lemma}\label{lm:glue} 
Let the situation be as above, and $\ball\subset S$ be a standard (analytic) ball. 
We consider the (analytic) open subset $\cU:=\rho(\cY_\ball)\subset X$. 
Then the following statements hold: 
\begin{enumerate}
\item There is a pointwise injective homomorphism 
$
{\bigl(\kern-.3ex
\underset{\mu\in M}{\mbox{$\bigoplus$}}\,
\eL_\mu^{\oplus m_\mu}
\bigr)}
{\otimes}\,\eO_\cU
\,{\to}\,\crl V{\otimes}\,\eO_\cU
$ 
whose restriction to $Y_s$ is the natural evaluation \eqref{eq:evm}, for all $s\in\ball$. 

\item $\crl V\otimes\eO_\cU$ is obtained as a successive extension of 
$\{\eL_j\}_{j\in J}\subset\Pic(X)$. 
\end{enumerate}
\end{m-lemma}

\begin{proof} \nit(i) 
First we prove that $\cU$ is indeed open. 
Since $\rho$ is an algebraic morphism, $\rho(\cY_\ball)\subset X$ is a locally closed 
analytic subvariety (it is constructible). If it is not open, the (local) components of 
$\rho|_{\cY_\ball}$ satisfy a non-trivial algebraic relation. This relation holds on 
whole $\cY$, since $\cY_\ball\subset\cY$ is open, which contradicts the hypothesis that 
$\rho(\cY)\subset X$ is open. 

Now we proceed with the proof of the lemma. For all $s\in\ball$, the restriction of  
$$
\ev:\kern-.5ex
\underset{\mu\in M}{\mbox{$\bigoplus$}}\,
\rho^*\eL_\mu\otimes\pi^*\pi_*\rho^*(\eL_\mu^{-1}\otimes\crl V)_\ball
\to (\rho^*\crl V)_\ball
$$
to $Y_s$ is the homomorphism \eqref{eq:evm}. 
The maximality of $\mu\in M$ implies that 
$$
\pi_*\rho^*(\eL_\mu^{-1}\otimes\crl V)\cong\eO_\ball^{\oplus m_\mu}, 
\forall\,\mu\in M,
$$
and $\ev$ is pointwise injective. We claim that, after suitable choices of bases in 
$\pi_*\rho^*(\eL_\mu^{-1}\otimes\crl V)$, the homomorphism $\ev$ descends to $\cU$. 
We deal with each $\mu\in M$ separately, the overall basis being the direct sum 
of the individual ones. 

Consider $\mu\in M$, and a base point $o\in\ball$. 
Then $\crl V':=\eL_\mu^{-1}\otimes\crl V$ has the properties:
\smallskip

-- $(\rho^*\crl V')_\ball\cong\eO_{\cY_{\ball}}^{\oplus m}\oplus
\underset{j\in J\sm\{\mu\}}{\bigoplus}
\rho^*{(\eL_\mu^{-1}\eL_j)}^{\oplus m_j}_\ball$. 

-- $\pi_*(\rho^*\crl V')_\ball\cong\eO_{\ball}^{\oplus m}$. 
We choose an isomorphism $\alpha_\ball$ between them.
\smallskip

-- $\pi^*\pi_*(\rho^*\crl V')_\ball\to (\rho^*\crl V')_\ball$ 
is pointwise injective; let $\eT\subset (\rho^*\crl V')_\ball$ be its image.
\medskip

\nit We choose a complement 
$\eW\,{\cong}\kern-.5ex\underset{j\in J\sm\{\mu\}}{\bigoplus}\kern-.5ex
\rho^*{(\eL_\mu^{-1}\eL_j)}^{\oplus m_j}$ of $\eT$, so 
${(\rho^*\crl V')}_\ball\,{=}\,\eT\,{\oplus}\,\eW.$ 
($\eW$ is defined up to $\Hom(\rho^*\crl V/\eT,\eT)$.) 
Then $\alpha_\ball$ above determines the pointwise injective homomorphism 
$\alpha\,{:}\,\eO_{\cY_{\ball}}^{\oplus m}\,{\to}\,(\rho^*\crl V')_\ball\,
{=}\,\eT{\oplus}\eW$ 
whose second component vanishes, since $\Gamma(\cY_{\ball},\eW)=0$. 
The left inverse $\beta\,{:}\,(\rho^*\crl V')_\ball\,{\to}\,\eO_{\cY_{\ball}}^{\oplus m}$ 
of $\alpha$ with respect to the splitting, satisfies $\alpha\circ\beta|_\eT=\bone_\eT$. 

\nit{\it Claim}\quad After a suitable change of coordinates 
in $\eO_{\cY_{\ball}}^{\oplus m}$, the homomorphisms $\alpha$ descends to 
$\rho(\cY_{\ball})\subset X$. Indeed, for any $s\in\ball$, we consider the diagram 
(recall that $Y_{os}\neq\emptyset$):
$$
\xymatrix@C=3em@R=1.25em{
\kern1ex\eO_{Y_{os}}^{\oplus m}\ar[r]^-{\alpha_o}\ar@{..>}[d]_-{a_s}
&
\crl V'_{Y_{os}}\ar@{=}[d]
&
\\ 
\kern1ex\eO_{Y_{os}}^{\oplus m}\ar[r]^-{\alpha_s}
&
\crl V'_{Y_{os}}\ar@<-2pt>`d[l] `[l]|-{\beta_s}
&
\text{with}\;
a_s:=\beta_s\circ\alpha_o\in\End(\mbb C^m).
}
$$
Similarly, we let $a'_s:=\beta_o\circ\alpha_s$. Then holds 
$a'_sa_s=\beta_o\alpha_s\beta_s\alpha_o=\beta_o\alpha_o=\bone$ 
(the second equality uses 
$\Img\bigl({\alpha_o|}_{Y_{os}}\bigr)
=\eT_{Y_{os}}=\Img\bigl({\alpha_s|}_{Y_{os}}\bigr)$, 
and $\alpha_s\beta_s|_\eT=\bone$), and similarly $a_sa'_s=\bone$. 
Thus $a_s\in\GL(m;\mbb C)$ for all $s\in\ball$, and the new trivialization 
$\tilde\alpha:=\alpha\circ a$ of $\eT$ satisfies  
\begin{equation}\label{eq:os}
\tilde\alpha_s=\tilde\alpha_o\text{ along }Y_{os},
\quad\forall\,s\in\ball, 
\end{equation}
because 
$\tilde\alpha_s|_{Y_{os}}{=}(\alpha_s\beta_s)\alpha_o|_{Y_{os}}
{=}\alpha_o|_{Y_{os}}{=}\tilde\alpha_o|_{Y_{os}}$. 

Also, for all $s,t\in\ball$ such that $Y_{st}\neq\emptyset$, 
the trivializations of $\eT_{Y_{st}}$ induced by $\tilde{\alpha}$ from 
$Y_s$ and $Y_t$, coincide; equivalently, the following diagram commutes: 
\begin{equation}\label{eq:st}
\xymatrix@R=1.5em@C=2.5em{
\eO_{Y_{st}}^{\oplus m}\ar[r]^-{\tilde\alpha_s}_-\cong\ar@{=}[d]
&
\eT_{Y_{st}}\ar@<-.5ex>@{=}[d]
\kern-1.5ex\ar@{}[r]|-{\mbox{$\subset$}}
&
\kern-1.5ex
\crl V'_{Y_{st}}\ar@<-2ex>@{=}[d]
&
\ar@{}[d]^-{\mbox{$
\Leftrightarrow\;
{\tilde\alpha_t^{-1}\circ\tilde\alpha_s|}_{Y_{st}}=\bone\in\GL(r;\mbb C).
$}}
\\ 
\eO_{Y_{st}}^{\oplus m}\ar[r]^-{\tilde\alpha_t}_-\cong
&
\eT_{Y_{st}}
\kern-1.5ex\ar@{}[r]|-{\mbox{$\subset$}}
&
\kern-1.5ex
\crl V'_{Y_{st}}
&
}
\end{equation}
Indeed, $Y_{ost}$ is non-empty and connected by \notr\kern-3pt, so is enough to prove 
that the restriction of \eqref{eq:st} to $Y_{ost}$  is the identity. 
After restricting \eqref{eq:os} to $Y_{ost}$, we deduce 
$$
\tilde{\alpha}_s|_{Y_{ost}}
=\tilde{\alpha}_o|_{Y_{ost}}
=\tilde{\alpha}_t|_{Y_{ost}}
\quad\Rightarrow\quad
{\tilde{\alpha}_t^{-1}\circ\tilde{\alpha}_s|}_{Y_{ost}}=\bone.
$$
Now we can conclude that the trivialization $\tilde\alpha$ of 
$\pi_*\rho^*(\eL_\mu^{-1}\otimes\crl V)$ descends to $\cU:=\rho(\cY_{\ball})$, 
as announced. Indeed, define 
\\[1ex] \centerline{
$\bar\alpha:\eO_\cU^{\oplus m}\to\crl V\otimes\eO_\cU$, 
$\bar\alpha(x):=\tilde\alpha_s(x)$ for some $s\in\ball$ such that  $x\in Y_s$. 
}\\[1ex]
The diagram \eqref{eq:st} implies that $\bar\alpha(x)$ is independent of 
$s\in\ball$ with $s(x)=0$. 

\nit(ii) Apply repeatedly the first part. 
\end{proof}

\begin{m-lemma}\label{lm:ext}
Let $Y_o\subset X$ be a subvariety and $\cU\supset Y_o$ be an open neighbourhood 
of it. Let $\eA,\eB,\eC$ be vector bundles on $X$, whose restriction to $\cU$ fit into 
$0\to\eA_\cU\to\eB_\cU\to\eC_\cU\to 0$. Assume that either one of the following 
conditions is satisfied:

\begin{tabular}{rrl}
&{\rm(i)} & $Y_o\subset X$ is $2^\apos$ (cf. \ref{def:ap>0});
\\ 
{\rm or}\hspace{1ex}&{\rm(ii)} & $\cd(X\sm Y_o)\les\dim X-3$. 
\end{tabular}

\nit Then $\eB$ is an extension of $\eC$ by $\eA$ on $X$. 
\end{m-lemma}

\begin{proof}
(i) As $Y_o$ is $2^\apos$, 
$\Ext^1(\eC,\eA)\to \Ext^1(\eC_{\tld Y_{m}},\eA_{\tld Y_{m}})$ is an isomorphism 
for an increasing sequence of thickenings $\{\tld Y_m\}_{m\ges\,m_0}$ of $Y_o$. 
The restriction of $\eB_\cU\in\Ext^1(\eC_\cU,\eA_\cU)$ to $\tld Y_m$ yields 
the extension $0\to\eA\to\eB'\to\eC\to 0$ over $X$, with the property that 
$\eB'_{\tld Y_m}\cong\eB_{\tld Y_m}$, for all $m\ges m_0$. 
As $Y_o$ is $2^\apos$, it follows $\eB\cong\eB'$. 

\nit(ii) The hypothesis that $\cd(X\sm Y_o)\les\dim X-3$ implies 
(cf. \cite[Ch. III, Theorem 3.4]{hart-as}) that for any locally free sheaf  
$\eF$ on $X$ the restriction below is an isomorphism: 
$$
H^i(X,\eF)\to H^i(\hat X_{Y_o},\hat{\eF}),
\;\text{for}\;i=0,1.
$$ 
Now use that 
$\Ext^1\big(\eC_{\hat X_{Y_o}},\eA_{\hat X_{Y_o}}\big)\cong 
H^1\big(\cHom(\eC_{\hat X_{Y_o}},\eA_{\hat X_{Y_o}})\big)$ 
and proceed as before. 
\end{proof}


\subsection{The splitting criterion} 
Let $\bbF\hra\bk$ be a finitely generated extension of $\mbb Q$, such that $X,\cY,\crl V$ 
are defined over $\bbF$; its algebraic closure $\bar\bbF\subset\bk$ is countable. 

\begin{m-theorem}\label{thm:split}
Let $X$ be a a smooth, projective variety, and assume the following: 
\begin{enumerate}
\item the situation is as in \ysx\hskip-1ex; 
\item 
$\crl V_{\mbb Y}$ splits on $\mbb Y$; 
alternatively, $\crl V_{Y_s}$ splits, for a very general $s\in S$;
\item 
$\mbb Y\subset X\otimes_{\bk}\bar\bbk=:X_{\bar{\Bbbk}}$ 
is either $2^\apos$ (e.g. it is $2^\pos$) or 
$\cd\big(X_{\bar\Bbbk}\sm\mbb Y\big)\les\dim X-3$.
\end{enumerate}
Then $\crl V$ is obtained by successive extensions of line bundles on $X$. 
If, moreover, $X$ has the property that $H^1(X,\eL)=0$ for all $\eL\in\Pic(X)$, 
then $\crl V$ splits. 

The very same statements remain valid if $X$ is defined over an uncountable, 
algebraically closed field, rather than over $\bbC$.
\end{m-theorem}

\begin{proof} 
First assume that $\crl V_{\mbb Y}$ splits. By lemma \ref{lm:ball}, there is a ball 
$\ball\subset S$, such that $(\rho^*\crl V)_{\cY_\ball}$ splits; lemma \ref{lm:glue} 
implies that $\crl V$ is a successive extension of line bundles on a tubular neighbourhood 
of $Y_o$, $o\in\ball$. It remains to apply lemma \ref{lm:ext}. 

Let $\tau\,{:}\,S\to S_{\bar\bbF}$ be the trace morphism; for $s\in S$, let 
$\Bbbk_0:=\bar\bbF\big(\tau(s)\big)$ be the residue field of $\tau(s)\in S_{\bar\bbF}$. 
For $s$ very general, $\tau(s)$ is the generic point of $S_{\bar\bbF}$, 
so $\Bbbk_0=\bar\bbF\big( S_{\bar\bbF}\big)$. 
Assume $\crl V_s$ splits; in the Cartesian diagram below 
$\crl V_{s}=\crl V_{\bar{\Bbbk}_0}\otimes\bbC$: 
$$
\xymatrix@R=1.5em@C=3em{
Y_s\ar[r]\ar[d]
&
\bbY_{\bar{\Bbbk}_0}
\ar[d]\ar[r]
&
\cY_{\bar\bbF}\ar[d]
\\  
\Spec(\bk)\ar[r]
&
\Spec(\bar\Bbbk_0)\ar[r]
&
S_{\bar\bbF}.
}
$$ 
The splitting of a vector bundle commutes with base change, for varieties defined over 
algebraically closed fields (cf. \cite{hal}). The previous discussion implies that 
$\crl V_{\bar{\Bbbk}_0}$ splits on $\bbY_{\bar{\Bbbk}_0}$ (the geometric generic 
fibre of $\cY_{\bar\bbF}\to S_{\bar\bbF}$); hence the same holds for 
$\crl V_{\bbY}$ on $\bbY=\bbY_{\bar\Bbbk_0}\otimes\bbC$. 
This brings us back to the previous case. For the final statement, 
we use once more that the splitting property commutes with base change. 
\end{proof}

The proof of the theorem even precises the meaning of the term `very general': 
if $\bbF$ is the field of definition of $X,\cY,\crl V$, then, in local affine coordinates 
coming from $S_{\bbF}$, the coordinates of $s\in S$ should be algebraically independent 
over $\bbF$. 

\begin{m-definition}\label{def:1split}
We say that the variety $X$ is \emph{$1$-splitting} if $H^1(X,\eL)=0$, 
for all $\eL\in Pic(X)$. 
\end{m-definition}
The simplest examples of $1$-splitting varieties are the Fano varieties of dimension at 
least two with cyclic Picard group, and products of such. 
In \ref{prop:homog-1split} are obtained examples of homogeneous $1$-splitting varieties 
with Picard groups of higher rank.

\begin{m-remark}\label{rmk} 
The genericity assumption in \ref{thm:split} plays an \emph{essential} role. 
If one is interested in the same result for arbitrary $Y$, one needs stronger 
positivity hypotheses, in order to apply effective cohomology vanishing results; 
in \cite{hal} I obtained a similar result for \emph{ample}, globally generated 
vector bundles. In this case, the situation \ysx holds automatically. 
\end{m-remark}
Theorem \ref{thm:split} will be illustrated in the following sections 
with concrete examples.  


\section{Positivity properties of zero loci of sections in vector bundles}\label{sct:0vb}

Throughout this section, $\eN$ stands for a \emph{globally generated} vector bundle 
on $X$ of rank $\nu$. 

\subsection{Sommese's $q$-ampleness for globally generated vector bundles}
\label{ssct:q-vb}

We briefly review the $q$-ampleness concept introduced in \cite{so}. 

\begin{m-proposition}{\rm (cf. \cite[Proposition 1.7]{so}).}\label{prop:q2}
The following statements are equivalent: 
\begin{enumerate}
\item The vector bundle $\eN$ is $q$-ample (cf. definition \ref{def:q-line}), 
that is for all coherent sheaves $\eF$ on $X$ holds 
\begin{equation}\label{eq:q21}
H^t(X,\eF\otimes\Sym^m(\eN))=0,\;\forall\,t\ges q+1,\;\forall\,m\gg0;
\end{equation}
\item The fibres of the morphism $\mbb P(\eN^\vee)\to |\eO_{\mbb P(\eN^\vee)}(1)|$ 
are at most $q$-dimensional. 
\end{enumerate}
\end{m-proposition}

\begin{m-proposition} 
\label{prop:q21} 
Assume $\eN$ is $q$-ample. 
Then $\eO_{\mbb P(\eN^\vee)}(1)$ is $q$-positive (as in \ref{def:p-line-bdl}).  
If $Y$ is the zero locus of a \emph{regular} section in $\eN$,  
then $Y\subset X$ is $(\dim X-\nu-q)^\pos$. 
\end{m-proposition}

\begin{proof}
The first statement is proved in \cite[Theorem 1.4]{mat}. 
Since the section is regular, $\codim_X(Y)=\nu$. 
For any coherent sheaf $\eF$ on $X$ holds 
$$
H^{\nu+t}(\mbb P(\eN),\eO(-\nu-m)\otimes\pi^*(\eF\otimes\det(\eN)^{-1}))
= H^{t+1}(X,\eF\otimes\Sym^m(\eN)),
$$
so the $H^{\ges \nu+q}$-cohomology on $\mbb P(\eN)$ vanishes; 
hence the same holds for $\tld X\subset\mbb P(\eN)$. 
\end{proof}

For $\nu=1$, the line bundle $\eN$ is $q$-ample if an only if the morphism 
$X\srel{}{\to}|\eN|$ has at most $q$-dimensional fibres. 
This property is easy to check, and convenient for concrete applications. 
In contrast, for $\nu\ges 2$, the criterion is not effective; the $q$-ampleness test 
for $\eN$ is too restrictive to check the positivity of the zero loci of its sections 
(cf. remarks \ref{rmk:sommese-weak}, \ref{rmk:sp-bad}.)


\subsection{The positivity criterion \ref{prop:x-b}}\label{ssct:x-b}

Suppose $Y\subset X$ is lci of codimension $\delta$, and the zero locus of a section 
$s$ in $\eN$, of rank $\nu\ges 2$; we \emph{do not assume} that $s$ is regular, 
so we allow $\delta<\nu$. 
In this context, the situation \ref{prop:x-b} arises as follows: 
since $Y$ is the zero locus of $s\in\Gamma(\eN)$, the blow-up $\tld X$ fits into 
\begin{equation}\label{eq:tld-x}
\xymatrix@R=1.5em@C=1.5em{
\tld X\ar@{^(->}[r]\ar[d]_-\pi\ar@<-5pt>[rrd]_-b
&
\mbb P(\eN)
{=}\,
\mbb P\Bigl(
\mbox{$\overset{\nu-1}{\bigwedge}$}\eN^\vee\otimes\det(\eN)
\Bigr)\ar@{^(->}[r]
&
X\times\mbb P\Bigl(
\mbox{$\overset{\nu-1}{\bigwedge}$}\Gamma(\eN)^\vee
\Bigr)\ar[d]
\\ 
X&&\mbb P:=\mbb P\Bigl(
\mbox{$\overset{\nu-1}{\bigwedge}$}\Gamma(\eN)^\vee
\Bigr),
}
\end{equation}
and holds 
\begin{equation}\label{eq:o1}
\eO_{\tld X}(E_Y)=\eO_{\mbb P(\eN)}(-1)\big|_{\tld X}=
\bigl(\det(\eN)\boxtimes\eO_{\mbb P}(-1)\bigr)\big|_{\tld X}.
\end{equation}

\begin{m-proposition}\label{prop:p}
Suppose $\det(\eN)$ is ample. If the dimension of the generic fibre of $b$ 
(over its image) is $p+1$, then $\eO_{\tld X}(E_Y)$ is $\dim b(\tld X)$-positive, 
and $Y$ is $p^\pos$. 
\end{m-proposition}

\begin{proof}
The assumptions of the proposition \ref{prop:x-b} are satisfied. 
\end{proof}

Observe that the propositions \ref{prop:q2} and \ref{prop:p} deal with 
complementary situations: $\eO_{\mbb P(\eN^\vee)}(1)$ is the pull-back of 
an ample line bundle; $\eO_{\tld X}(E_Y)$ is relatively ample for some morphism.

Let $W\subseteq\Gamma(\eN)$ be a vector subspace which generates $\eN$; 
let $\dim W=\nu+u+1$. 
A globally generated vector bundle $\eN$ on $X$ is equivalent to a morphism 
$f{:}\,X{\to}\,\Grs(W;\nu)$ into the Grassmannian of $\nu$-dimensional quotients 
of $W$; $\det(\eN)$ is ample if and only if $\vphi$ is finite onto its image. 

For $\Grs(W;\nu)$ and $\eN$ the universal quotient bundle on it, we can explicitly 
write the morphism $b$ in \eqref{eq:tld-x}: $\mbb P(\eN)\to\mbb P$ is defined by 
\begin{equation}\label{eq:q}
\mbox{$
(x,\lran{e_x})\mt 
\det(\eN_x/\lran{e_x})^\vee\subset
\overset{\nu-1}{\bigwedge}\eN_x^\vee\subset
\overset{\nu-1}{\bigwedge}W^\vee.
$}
\end{equation}
($\lran{e_x}$ stands for the line generated by $e_x\in\eN_x$, $x\in\Grs(W;\nu)$.) 
The restriction to the Grassmannian corresponds to the commutative diagram  
\begin{equation}\label{eq:wn}
\xymatrix@R=1.75em{
0\ar[r]&\eO_{\Grs(W;\nu)}\ar[r]^-{s}\ar@{=}[d]&
W\otimes\eO_{\Grs(W;\nu)}\ar[r]\ar@{->>}[d]^-{\;\beta}&
W/\lran{s}\otimes\eO_{\Grs(W;\nu)}\ar[r]\ar@{->>}[d]&0
\\ 
&\eO_{\Grs(W;\nu)}\ar[r]^-{\beta s}&\eN\ar[r]&\eN/\lran{\beta s}\ar[r]&0.
}
\end{equation}
Thus $b$ is the desingularization of the rational map 
\begin{equation}\label{eq:q-grs}
g_s:\Grs(W;\nu)\dashto\Grs(W/\lran{s};\nu-1),\quad 
[W\surj N]\mt [W/\lran{s}\;\surj N/\lran{\beta s}]
\end{equation}
followed by the Pl\"ucker embedding of $\Grs(W/\lran{s};\nu-1)$; 
the indeterminacy locus of $b$ is $\Grs(W/\lran{s};\nu)\subset\Grs(W;\nu)$.

\begin{m-remark}\label{rmk:sommese-weak}
We mentioned that Sommese's $q$-ampleness criterion is not effective for $\nu\ges 2$. 
For $X=\Grs(\nu+u+1;\nu)$, Sommese's criterion implies that $\eN$ is $q$-ample, 
with $q=\dim\mbb P(\eN^\vee)-\mbb P^{\nu+u}=\dim X-(u+1)$; 
hence $Y=\Grs(\nu+u;\nu)$, the zero locus of a generic section of $\eN$, is $(u+1-\nu)^\pos$. 
On the other hand, the criterion \ref{prop:x-b} implies that $Y$ is actually $u^\pos$. 
\end{m-remark}

There is a `universal' rational map $g_{\rm univ}$ containing the maps $g_s$ above, 
as $s$ varies: 
\begin{equation}\label{eq:G-grs}
\begin{array}{c}
g_{\rm univ}:\mbb P(W)\times\Grs(W;\nu)\dashto\Flag(W;\nu+u,\nu-1),
\\[1ex] 
(\lran{s},[W\surj N])\mt 
\bigl[W\surj\frac{W}{\lran{s}}\surj\frac{N}{\lran{\beta(s)}}\bigr]. 
\end{array}
\end{equation}
(The right hand side denotes the flag variety of successive quotients of $W$.)
It is undefined on the incidence variety 
$\cal J:=\{\,(\lran{s},N)\mid s\in\Ker(W\surj N)\}$. 

Back to the general case of a globally generated vector bundle on a variety $X$. 
By varying $s\in W$, one obtains the family of subvarieties of $X$ (over $\mbb P(W)$) 
$$
\cY:=(\mbb P(W)\times X)\times_{\mbb P(W)\times\Grs(W;\nu)}\cal J,
$$ 
and the situation mentioned at \ref{rmk:x-b}:
\begin{equation}\label{eq:xg}
\xymatrix@R=2em@C=3em{
\Bl_\cY(\mbb P(W)\times X)
\ar@{->}[r]\ar[d]_-\pi\ar@<-3pt>[drr]_<(.2)b|!{[r];[dr]}\hole
&
\Bl_{\cal J}\bigl(\mbb P(W)\times\Grs(W;\nu)\bigr)\ar[dr]^-{\tld g}\ar[d]
&
\\ 
\mbb P(W)\times X\ar@{->}[r]^-{\vphi}
&
\mbb P(W)\times\Grs(W;\nu)\ar@{-->}[r]_-{g_{\rm univ}}
&
\Flag(W;\nu+u;\nu-1).
}
\end{equation}
(The intersection $X\cap\Grs(W/\lran{s};\nu)$ may not be transverse 
for certain $s\in W$.) Proposition \ref{prop:p} can be restated as follows.

\begin{m-corollary}\label{cor:grs}
Let $\vphi:X\to\Grs(W;\nu)$, $\nu\ges 2$, be a morphism finite onto its image. 
If the general fibre of 
$g_{\rm univ}\circ\vphi:\mbb P(W)\times X\dashto\Flag(W;\nu+u;\nu-1)$ 
is at least $(p+1)$-dimensional, then $Y_s$ is $p^\pos$, for all $s\in W$ 
such that $Y_s$ is irreducible lci in $X$. 
\end{m-corollary}


\subsection{The Picard group and the diagram \ysx}\label{ssct:diag-0loci} 

The sheaf $\euf K$ defined by 
$$
0\to\euf K:=\Ker(\ev)\to W\otimes\eO_X\srel{\ev}{\to} \eN\to 0
$$
is locally free, and the incidence variety 
\begin{equation}\label{eq:yn}
\cY=\{(\lran{s},x)\mid s(x)=0\}\subset X\times\mbb P(W)
\end{equation} 
is isomorphic to $\mbb P(\euf K)$. We denote $\pi$ and $\rho$ the projections onto 
$\mbb P(W)$ and $X$ respectively; actually we restrict ourselves to a sufficiently 
small open subset $S\subset\mbb P(W)$.  
The vector bundles $\eN^{\oplus2}$ and $\eN^{\oplus3}$ are generated by 
$W^{\oplus2}$ and $W^{\oplus3}$ respectively, and determine the double and triple 
self-intersection diagrams \eqref{eq:y23}. 

\begin{m-proposition}\label{prop:pic-0loci}
Let $\cY\subset S\times X$ be as in \eqref{eq:yn}, with  $S\subset\mbb P(W)$ 
suitably small. Then the conditions of \ysx are fulfilled as soon as $\eN$ 
satisfies any of the following conditions: 
\begin{enumerate}
\item 
$\eN$ is Sommese-$q$-ample (cf. \ref{prop:q2}), and \newline 
\text{$\dim X-q\ges 3\nu+1$ (for $\nu\ges 2$), or $\dim X-q\ges 5$ (for $\nu=1$).}

\item[or\phantom{l}]
\item 
$\nu\ges 2$, and the generic fibre of $\,\mbb P(W)\times X\dashto\Flag(W;\nu+u;\nu-1)$ 
is at least $2(\nu+1)$-dimensional (cf. \eqref{eq:xg}). 
\end{enumerate} 
\end{m-proposition}

\begin{proof} 
(i) The non-emptiness of the triple intersections requires $\dim X-q\ges 3\nu+1$,  
by lemma \ref{lm:yy'}; 
the isomorphism of Picard groups requires $\dim X-q\ges 2\nu+3$, 
by theorem \ref{thm:pic}. 

\nit(ii) In this case, $Y_s\subset X$ is $(2\nu+1)^\pos$, so the triple intersections 
are non-empty (cf. \ref{lm:yy'}); 
the isomorphism of Picard groups requires $(2\nu+1)-\nu\ges 3$ (cf. \ref{lm:ball}).
\end{proof}


\section{Splitting along zero loci of globally generated vector bundles}
\label{sct:split-0loci}

Let $\crl V$ be an arbitrary vector bundle on $X$. 
The previous discussion immediately yields the following criterion. 

\begin{m-theorem}\label{thm:split-vb1}
Let $X$ be $1$-splitting variety (cf. \ref{def:1split}), $\eN$ be a globally generated 
vector bundle of rank $\nu$ on $X$ such that $\det(\eN)$ is ample, and 
$W\subset\Gamma(X,\eN)$ be a generating vector subspace. 
Assume that $\eN$ satisfies one of the conditions in proposition \ref{prop:pic-0loci}. 

Then $\crl V$ splits on $X$, if and only if its restriction to the zero locus 
of a very general $s\in W$ splits. 
(If $X$ is not $1$-splitting, $\crl V$ is a successive extension of line bundles.)
\end{m-theorem}

\begin{m-example}\label{expl:xgr}
Consider $\Grs(W;\nu)$, with $\nu\ges 2$, $\dim W=\nu+u+1$, and let 
$X\subset\Grs(W;\nu)$ be an arbitrary subvariety (e.g. the Grassmannian itself) 
such that 
$$
\Pic(X)\cong\mbb Z,\quad\codim_{\Grs(W;\nu)}(X)\les u-(2\nu+1).
$$
Then a vector bundle $\crl V$ splits if and only if it does so along 
a very general zero locus $Y\subset X$ of a section in $\eN_X$. 
Indeed, in this case proposition \ref{prop:pic-0loci}(ii) is satisfied. 
Moreover, observe that $\Pic(Y)\cong\Pic(X)\cong\mbb Z$ so the procedure can be 
iterated as long as the codimension condition above is satisfied. 
\end{m-example}
In section \ref{sct:grass}, this principle will be further illustrated with numerous 
concrete examples. 

Is natural to ask what happens if one drops the hypothesis that $\det(\eN)$ is ample. 
The Stein factorization of $\vphi:X\to\Grs(W;\nu)$ decomposes it into a morphism 
with connected fibres followed by a finite map. 

\begin{m-theorem}\label{thm:split-vb2} 
Let $\vphi:X\to X'$ be a smooth morphism of relative dimension $d\ges 1$, 
$\eN'$ a globally generated vector bundle on $X'$ of rank $\nu$ 
with $\det(\eN')$ ample, and $\eN=\vphi^*\eN'$. 

Assume that $X$ is $1$-splitting, and $\eN'$ satisfies proposition \ref{prop:pic-0loci}. 
Then  $\crl V$ splits on $X$ as soon as $\crl V$ splits along  $Y_s$, for $s$ very general.
\end{m-theorem}

\begin{proof}
The conditions of theorem \ref{thm:split} are fulfilled. 
Indeed, consider the family $\cY'$ of subvarieties of $X'$ as above 
and let $\cY:=\vphi^{-1}(\cY')$. The positivity is preserved under pull-back 
(cf. lemma \ref{lm:pull-back}, \ref{lm:pull-back-smooth}). 
For the condition (Pic), the morphism $\tld\vphi:\tld X\to\tld X'$ induced 
at the level of the blows-up is still smooth of relative dimension $d$, 
so $\eO_{\tld X}(E_Y)=\tld\vphi^*\eO_{\tld X'}(E_{Y'})$ is $(q+d)$-ample; 
it remains to apply \ref{thm:pic}. 
\end{proof}

It is surprising that this yields new results even in the simplest case $X=X'\times V$, 
with $X',V$ smooth, and $\eN'=\eO_{X'}(1)$ is globally generated, ample. 
If $X$ is $1$-splitting, then $X',V$ are both $1$-splitting; if either $X'$ or $V$ 
are simply connected, the converse is true.

\begin{m-corollary}\label{cor:restr-div}
Assume that $X=X'\times V$ is $1$-splitting, $\dim X'\ges 5$, 
and $\eO_{X'}(1)$ is ample, globally generated. Then $\crl V$ splits on $X$, 
if it splits along a very general divisor in $|\vphi^*\eO_{X'}(1)|$. 
\end{m-corollary}


\section{Positivity properties of sources of $G_m$-actions}\label{sct:fixed} 

Another (totally different) framework which leads to the situation \ysx arises 
in the context of the actions of the multiplicative group on (almost) homogeneous 
varieties. 

\subsection{Basic properties of the BB-decomposition}\label{ssct:bb}

We start with general considerations which should justify the appearance of the 
homogeneous varieties in the next section. 
Let $G$ be a connected reductive linear algebraic group, and $T\subset B\subset G$ 
be a maximal torus and a Borel subgroup. 
Finally, let $X$ be a smooth projective $G$-variety with a faithful $G$-action 
$\mu:G\times X\to G\times X$. 

Also, consider a $1$-parameter subgroup ($1$-PS for short) $\l:G_m\to G$, 
where $G_m=\bk^*$ is the multiplicative group; we assume $\l(G_m)\subset T$. 
The $1$-PS induces the action $\l:G_m\times X\to X$, which determines the so-called 
Bialynicki-Birula (BB for short) decomposition of $X$ into locally closed subsets. 
Below are summarized its basic properties (cf. \cite{bb,kon}): 
\begin{itemize}
	\item 
		The specializations at $\{0,\infty\}=\mbb P^1\sm G_m$ are denoted 
		$\underset{t\to 0}{\lim}\l(t)\times x$ and $\underset{t\to\infty}{\lim}\l(t)\times x$; 
		they are both fixed by $\l$. 
	\item 
		The fixed locus $X^\l$ of the action is a disjoint union 
		$\underset{s\in S_\BB}{\coprod}\kern-1exY_s$ of smooth subvarieties. 
		For $s\in S_\BB$, 
		$Y_s^+:=\{x\in X\mid\underset{t\to 0}{\lim}\l(t)\times x\in Y_s\}$ 
		is locally closed in $X$ (a BB-cell).
	\item 
		$X=\underset{s\in S_\BB}{\coprod}\kern-1exY_s^+$, and the morphism 
		$Y_s^+\to Y_s,\; x\mt\underset{t\to0}{\lim}\l(t)\times x$ 
		is a locally trivial, affine space fibration; it is not necessarily a vector bundle. 
	\item 
		The \emph{source} $Y:=Y_{\rm source}$ and 
		the \emph{sink} $Y_{\text{\rm sink}}$ of the action are uniquely characterized 
		by the fact that $Y^+=Y_{\rm source}^+\subset X$ is open 
		and $Y_{\rm sink}^+=Y_{\rm sink}$. 
	\item 
		By composing $\l$ with the involution $t\mt t^{-1}$ of $G_m$, one gets 
		the so-called \emph{minus} BB-decomposition 
		$\; X=\underset{s\in S_\BB}{\coprod}\kern-1exY_s^-.$
\end{itemize}
We denote: 
\begin{equation}\label{eq:l}
\begin{array}{rl}
G(\l)&:=\big\{g\in G\mid g^{-1}\l(t)g=\l(t),\;\forall\,t\in G_m\big\}
\text{ the centralizer of $\l$ in $G$}, 
\\[1ex] \disp
P(\pm\l)&:=\big\{g\in G\mid\lim_{t\to 0}
\big(\l(t)^{\pm}g\,\l(t)^{\mp}\big)\text{ exists in }G\big\},
\\[1ex] \disp 
U(\pm\l)&:=\big\{g\in G\mid\lim_{t\to 0}
\big(\l(t)^{\pm}g\,\l(t)^{\mp}\big)=e\in G\big\}.
\end{array}
\end{equation}
Then $G(\l)$ is a connected, reductive subgroup of $G$, $P(\pm\l)\subset G$ are 
parabolic subgroups,  $G(\l)$ is their Levi-component, and $U(\pm\l)$ the unipotent 
radical (cf. \cite[\S 13.4]{sp}). 

\begin{m-lemma}\label{lm:y-invar}
\nit{\rm(i)} $Y$ is invariant under $P(-\l)$, thus for $G(\l)$ too. 

\nit{\rm(ii)} $Y_s^+$ is $P(\l)$-invariant and $U(\l)$ preserves the fibration 
$Y_s^+\to Y_s$, for all $s\in S_\BB$. 
\end{m-lemma}

\begin{proof} 
See \cite{hal+taj}. 
\end{proof}

Theorem \ref{thm:split} requires that the embedded deformations of $Y$ sweep out 
an open subset of $X$; 
this is achieved if $GY$, defined set theoretically as $\{gy\mid g\in G,\;y\in Y\}$, 
is open in $X$. 

\begin{m-lemma}\label{lm:gy}
$GY$ has a natural structure of a closed subscheme of $X$. 
Therefore $GY\subset X$ is open if and only if $GY=X$. 
\end{m-lemma}

\begin{proof}
Indeed $GY$ is the image of $\mu:G\times Y\to X$. 
Since $Y$ is $P(-\l)$-invariant, it factorizes 
$(G\times Y)/P(-\l)\to X$, for $p\times(g,y):=(gp^{-1},py)$. 
But $P(-\l)$ is parabolic, so $(G\times Y)/P(-\l)$ is projective, 
hence $\Img(\mu)$ is a closed in $X$. 
\end{proof}

\begin{m-remark}\label{rmk:open}
Since $P(\l)P(-\l)Y=P(\l)Y=U(\l)G(\l)Y=U(\l)Y$, and  
$P(\l)P(-\l)$ is open in $G$, we deduce: 
\\ \centerline{
$GY=X\;\Leftrightarrow\;U(\l)Y\subset X$ is open. 
}\\[.5ex]
This observation hints towards the fact that the $G$-varieties satisfying the 
lemma \ref{lm:gy} should be homogeneous (or, at least, have an open $B$-orbit). 
\end{m-remark}


\subsection{The positivity of $Y\subset X$}\label{ssct:posYX}

We follow the same steps as in the section \ref{sct:0vb}: determine the positivity 
of the source of a $G_m$-action. In section \ref{sct:q-tech} we obtained two methods 
for doing this.  


\subsubsection{Apply the proposition \ref{prop:trans}}\label{sssct:trans}

Let us recall that $Y\subset X$ is $p^\pos$ as soon as:  

-- $\eN_{Y/X}$ is $(\dim X-p)$-ample; 

-- $\cd(X\sm Y)\les\dim X-(p+1)$.

\nit It does not seem to exist a uniform answer for the amplitude of the normal bundle; 
in contrast, there is a compact formula for the cohomological dimension. 

\begin{m-proposition}\label{prop:cdX-Y}
If $Y\subset X$ is the source of a $G_m$-action, then holds 
$$
\cd(X\sm Y)=\dim(X\sm Y^+).
$$
\end{m-proposition}\vspace{-10pt}
\begin{equation}\label{eq:xy+}\kern-1\parindent
\begin{array}{lll}
\text{Observe that:}\quad&
\dim X-\dim (X\sm Y^+)
&
=\dim X-\max\{\dim\ovl{\,Y_s^+}\,\big|\, s\neq{\rm source}\}
\\[1ex] 
&&
=
\underset{s\neq{\rm source}}{\min}
\Big\{
				\begin{array}{l}
					\text{number of strictly negative}\\ 
					\text{weights of $\l$ on $T_X|_{Y_s}$}
				\end{array}
\Big\}.
\end{array}
\end{equation}

\begin{proof} 
See \cite{hal+taj}. 
\end{proof}


\subsubsection{Apply the proposition \ref{prop:x-b}}\label{sssct:x-b}
We remark that the $G_m$-action leads to the situation analyzed 
in proposition \ref{prop:x-b}. By linearizing the $G_m$-action in a very ample 
line bundle on $X$, one gets a $G_m$-equivariant embedding of $X$ into some 
$\mbb P^N$, such that $X$ is not contained in a hyperplane. 
In coordinates $z_{N_0}\in\bk^{N_0},\ldots,z_{N_r}\in\bk^{N_r}$, 
the $G_m$-action on $\mbb P^N$ is: 
\begin{equation}\label{eq:c*}
t\times[z_{N_0},z_{N_1},\ldots,z_{N_r}]=
[z_{N_0},t^{m_1}z_{N_1},\ldots,t^{m_r}z_{N_r}],
\quad\text{with }0<m_1<\ldots<m_r. 
\end{equation}
The source and sink of $\mbb P^N, X$ are respectively: 
\begin{equation}\label{eq:YP}
\begin{array}{lcll}
\mbb P^N_{\rm source}=\{[z_{N_0},0,\ldots,0]\}, 
&&
\mbb P^N_{\rm sink}=\{[0,\ldots,0,z_{N_r}]\},
\\[1ex] 
Y=Y_{\rm source}=X\cap \mbb P^N_{\rm source},
&& 
Y_{\rm sink}=X\cap \mbb P^N_{\rm sink},
\\[1ex]
Y^+=X\cap (\mbb P^N_{\rm source})^+,
&&
(\mbb P^N_{\rm source})^+
=\{\uz=[z_{N_0},z_{N_1},\ldots,z_{N_r}]\mid z_{N_0}\neq 0\}.
\end{array}
\end{equation}
(The intersections are set theoretical.) 
Note that $\mbb P^N_{\rm source}$ is the indeterminacy locus of  the rational map 
$$
\mbb P^N\dashto\mbb P^{N'},\quad
[z_{N_0},z_{N_1},\ldots,z_{N_r}]\mt[z_{N_1},\ldots,z_{N_r}],
$$
which can be resolved by a simple blow-up. 
By restricting to $X$, we get the diagram:
\begin{equation}\label{eq:y-source}
\xymatrix@R=1.5em@C=3em{
\tld X=\Bl_{Y}(X)\ar[d]\ar[r]^-{\tld\iota}&
\wtld{\mbb P^N}:=\Bl_{\mbb P^N_{\rm source}}(\mbb P^N)\ar[d]\ar[r]^-b&
\mbb P^{N'}
\\ 
X\ar[r]^-\iota&\mbb P^N&&
}
\end{equation}
If $E\subset\wtld{\mbb P^N}$ is the exceptional divisor, 
$$
\eO_{\tld X}(E_Y)={\eO_{\wtld{\mbb  P^N}}(E)\big|}_{\tld X}
={\eO_b(1)\otimes b^*\eO_{\mbb P^{N'}}(-1)\big|}_{\tld X}\,,
$$ 
and $\eO_b(1)$ is relatively ample. 
Thus, the requirements of \ref{prop:x-b} are satisfied; to compute the positivity 
of $Y\subset X$, one must estimate the dimension of $(b\tld\iota)(\tld X)$. 
In general, it is not clear how to do this, but we obtain an appealing result under 
additional hypotheses.

\begin{m-proposition}\label{prop:y-sink}
Let $X$ be a smooth $G_m$-variety with source $Y$, and assume that 
$G_m$ acts on the normal bundle $\eN_{Y/X}$ by scalar multiplication. 
Then $Y\subset X$ is $p^\pos$, with 
\\[.5ex]\centerline{$p=\dim X-\dim(X\sm Y^+)-1.$} 
\end{m-proposition}

\begin{proof}
The first step is to show that $X\subset\mbb P^N$ is invariant under the $G_m$-action 
by \emph{scalar multiplication} on the coordinates $(z_{N_1},\ldots,z_{N_r})$. 
Indeed, since $G_m$ acts by scalar multiplication on $\eN_{Y/X}$, it follows that 
$Y^+\to Y$ is a vector bundle (cf. \cite[Remark pp. 491]{bb}), so $Y^+$ is the 
total space of $\eN_{Y/X}$; for clarity we denote it by 
$\unbar{\sf N}_{Y/X}:=\Spec\big(\Sym^\bullet\eN_{Y/X}^\vee\big)$. 
Thus $\unbar{\sf N}_{Y/X}$ is a Zariski open subset of $X$, on which 
$G_m$-acts by scalar multiplication. Moreover, the inclusions 
$$
\unbar{\sf N}_{Y/X}\subset\unbar{\sf N}_{\mbb P^N_{\rm source}/\mbb P^N}
\srel{\eqref{eq:YP}}{=}
\{\,[z_{N_0}:z_{N_1}:\ldots]\mid z_{N_0}\neq0\,\}\subset\mbb P^N
$$ 
are $G_m$-equivariant. But the diagonal multiplication on $(z_{N_1},\ldots,z_{N_r})$ 
exists on the whole $\mbb P^N$, so $X\subset\mbb P^n$ is $G_m$-invariant too. 

Now we proceed to estimate the dimension of $(b\tld\iota)(\tld X)$. 
The previous step implies that for all $\uz=[\uz_0:\uz']\in Y^+\sm Y$ holds 
$\disp[\,0:\uz'\,]=\lim_{t\to\infty}[\,\uz_0:t\cdot\uz'\,]\in X\sm Y^+.$ 
This implies that $b(\tld X)=b(X\sm Y^+)$, as desired. 
\end{proof}


\subsection{The Picard group and the diagram \ysx}\label{ssct:pic-source}

\begin{m-proposition}\label{prop:pic-source} 
Let the situation be as in proposition \ref{prop:y-sink}, with 
$\dim X-\dim(X\sm Y^+)\ges 2$. Then $\Pic(X)\to\Pic(Y)$ is an isomorphism. 
\end{m-proposition}

\begin{proof} 
The inclusion $Y^+\subset X$ yields the exact sequence 
\begin{equation}\label{eq:picY}
0\to\Pic(X\sm Y^+)
=\Pic\big(\underset{s\neq\;{\rm source}}{\bigcup}\kern-1.5ex\ovl{\,Y_s^+}\;\big)
\to\Pic(X)\to\Pic(Y^+)\cong\Pic(Y)\to 0.
\end{equation}
The isomorphism on the right hand side holds because $Y^+\to Y$ is an affine space 
fibration. 

The left hand side is the free abelian group generated by the divisors  
$\ovl{\,Y_s^+\!}$, $s\in S_\BB\sm\{\rm source\}$. 
But the hypothesis implies that all the components of $X\sm Y^+$ have codimension 
at least two in $X$, so $X\sm Y^+$ contains no divisors. 
\end{proof}

By using the $G$-action on $X$, we can deform $Y$ to $Y_g:=gY$, with $g\in G$; 
the latter is the source of the $G_m$-action for the 1-PS $\l_g:=\Ad_g(\l)$. 
If $G$ is sufficiently large (see the remark \ref{rmk:open}), 
the family $\{Y_g\}_{g\in G}$ sweeps out an open subset of $X$. 
In order to apply the theorem \ref{thm:split}, one still has to control 
the Picard group of the double intersections $Z_g:=Y\cap Y_g$. 
One way to achieve this is by proving sufficient positivity of $Y$ 
(cf. subsection \ref{sssct:x-b}) and applying theorem \ref{thm:pic}. 
However, if this method fails, one has to use additional symmetry. 

\begin{m-corollary}\label{cor:gm2}
Let $X$ be a $G$-variety and $\l$ a 1-PS of $G$; consider $\gamma\in\Weyl(G)/\Weyl(G(\l))$. We assume the following: 
\begin{enumerate}
\item 
The sources $Y,Y'{=}\,\gamma Y$ of $\l,\l'$ intersect transversally; 
denote $Z:=Y\cap Y'$; 
\item 
$\codim_X(Y)+2\les\dim X-\dim(X\sm Y^+)$.
\end{enumerate}
Then, for generic $g\in G$, holds: 
\begin{enumerate}
\item[--] 
The intersection $Z_g:=Y\cap Y_g$ is transverse.
\item[--] 
The restrictions $\Pic(X)\to\Pic(Y)\to\Pic(Z_g)$ are isomorphisms.
\end{enumerate}
\end{m-corollary}

\begin{proof}
It is enough to prove the claims for $g=\gamma$, because both are unchanged under 
small perturbations. The isomorphism $\Pic(X)\to\Pic(Y)$ has been proved before. 
For the isomorphism $\Pic(Y)\to\Pic(Z)$, we observe that $\l'$ leaves $Y$ invariant, 
so $Y$ admits itself a \BB-decomposition. Since $Y\sm Z=Y\cap(X\sm Z)$, 
we deduce that 
\\[.25ex]\centerline{$
\dim(Y\sm Z^+)=\cd(Y\sm Z)\les\cd(X\sm Z)=\cd(X\sm Y)
=\dim(X\sm Y^+)\srel{\rm(ii)}{\les}\dim(Y)-2.
$}\\[1ex]
Hence $Y\sm Z^+$ contains no divisors and the exact sequence \eqref{eq:picY} 
yields the conclusion. 
\end{proof}

\begin{m-theorem}\label{thm:pic-sources} 
{\rm(i)} 
Consider the $G$-subvariety $\cY:=\mu(G\times Y)\subset G\times X$. 
Assume that: 
\\[1ex] \centerline{
$\cY\srel{\rho}{\to}X\text{ is smooth,}
\quad\dim X-\dim(X\sm Y^+)\ges 2\codim_X(Y)+2\ges 6$.
}\\[1ex]
Then the conditions \ysx are satisfied. 

\nit{\rm(ii)} Assume furthermore that $Z_g:=Y\cap Y_g$ has the property 
that $\Pic(Y)\to\Pic(Z_g)$ is an isomorphism, for $g\in G$ generic. 
Then the splitting criterion \ref{thm:split} applies. 
\end{m-theorem}
Note that the condition (ii) can be settled in the situations described in 
\ref{prop:y-sink} and \ref{cor:gm2}.

\begin{proof} 
(i) Since $\rho$ is smooth, so are the generic self intersections of $\cY$. 
These double and triple self-intersections are also non-empty and connected 
(cf. \ref{lm:yy'}); thus the \arm and \notr conditions in \ysx are satisfied. 

\nit(ii) By proposition \ref{prop:cdX-Y}, we have 
$\cd(X\sm Y)=\dim(X\sm Y^+)<\dim X-3$.  
\end{proof}


\section{Splitting criteria for vector bundles on homogeneous varieties}
\label{sct:split-homog}

In this section we specialize the previous discussion to homogeneous varieties. 
Assume $X=G/P$, where $G$ is connected, reductive, and $P$ is a parabolic 
subgroup, and consider a $1$-PS $\l$ of $G$. For any parabolic subgroup 
$Q$ of $G$, denote $\Weyl(Q):=\Weyl\big(\Levi(Q)\big)$. 

The adjoint action of $\l$ on $\Lie(G)$ decomposes it into the direct sum of its weight 
spaces; we group them into the zero, strictly positive and negative weight spaces: 
\begin{equation}\label{eq:LG}
\Lie(G)=\Lie(G)^0_\l\oplus\Lie(G)^{+}_\l\oplus\Lie(G)^{-}_\l, 
\;\text{and}\;
\Lie(P(\pm\l))=\Lie(G)^0_\l\oplus\Lie(G)^\pm_\l. 
\end{equation}

\begin{m-lemma}\label{lm:y-homog} 
The following statements hold:
\\ \nit{\rm(i)} 
The components of $X^\l$ are homogeneous for the action of $G(\l)$.
\\ \nit{\rm(ii)} 
The source $Y$ contains $\hat e\in G/P$ if and only if 
$\l\subset P$ and $\Lie(G)^-_\l\subset\Lie(P)$. 
\end{m-lemma}

\begin{proof}
See \cite{hal+taj}. 
\end{proof}

The maximal torus $T$ acts with isolated fixed points on $G/P$; they are precisely 
$wP$, with $w\in\Weyl(G)\slash\Weyl(P)$. For any ${s\in S_\BB}$, the component 
$Y_s\subset X^\l$ is $G(\l)$-invariant and $T\subset G(\l)$, so $Y_s^T\neq\emptyset$; 
conversely, $\l$ fixes $(G/P)^T$, so $(G/P)^T=\underset{s\in S_\BB}{\bigcup}Y_s^T$. 

Now we recall some classical facts about Bruhat decompositions (cf. \cite[\S 8]{mitch}, 
\cite[\S 8]{sp}). 
If $Q,P\subset G$ are two parabolic subgroups, $G/P$ decomposes into locally closed $Q$-orbits:
$$
G/P=\underset{w\in S_\Bht}{\mbox{$\coprod$}}\kern-1exQwP,
\;\text{with}\;
S_\Bht=\Weyl(Q)\big\backslash\Weyl(G)\,\big\slash\Weyl(P).
$$
Actually, $S_\Bht$ parameterizes the $\Weyl(Q)$-orbits in $(G/P)^T$. 
Each double coset in $S_\Bht$ contains a unique representative of minimal length; 
for each $w\in S_\Bht$ of minimal length, 
$$
\dim (QwP)={\rm length}(w)+\dim\big(\Levi(Q)\big/\Levi(Q)\cap wPw^{-1}\big).
$$

\begin{m-proposition}\label{prop:BB-homog}
The Bialynicki-Birula decomposition of $G/P$ for the action of $\l$ coincides with 
the Bruhat decomposition for the action of $P(\l)$. 
\end{m-proposition}

\begin{proof}
See \cite{hal+taj}. 
\end{proof}

For homogeneous varieties, the criterion \ref{thm:split} yields the following.

\begin{m-theorem}\label{thm:split-homog} 
Let $X=G/P$ and $\l$ be a 1-PS of $T\subset G$; we denote by $Y$ its source. 
Consider also $\gamma\in\Weyl(G)/\Weyl(G(\l))$. We assume the following: 
\begin{enumerate}
\item 
$X$ is $1$-splitting; 
\item 
$Y, Y'=\gamma Y$ intersect transversally;
\item 
$\dim X-\dim(X\sm Y^+)\ges 2\codim_X(Y)+2\ges 6$.
\end{enumerate}
Then an arbitrary vector bundle $\crl V$ on $G/P$ splits if and only if 
$g^*\crl V$ splits on $Y$, for a very general $g\in G$. 
\end{m-theorem}
A pleasant feature is that the splitting of $\crl V$ is reduced to the splitting along 
a homogeneous subvariety of $G/P$, so the procedure can be iterated 
(cf. \ref{cor:iterate} below). Explicit calculations are performed 
in section \ref{sct:grass}. 

\begin{m-remark}\label{rmk:gamma}
The element $\gamma$ above is typically a simple reflection in 
$\Weyl(G)\sm\Weyl(G(\l))$. 
Consider, for instance, $X=\Grs(u;w)$, $G=\GL(w)$, 
and $\l(t)=\diag[t^{-1},1,\ldots,1]$. The source of the action is 
$Y=\{U\subset\mbb C^w\mid e_1=\lran{1,0,\ldots,0}\in U\}.$ 
The appropriate $\gamma$ equals 
$\diag\bigg[ \Big[\begin{array}{cc} 0,1\\1,0 \end{array}\Big],1,\ldots,1\bigg]$; 
it is the reflection corresponding to the simple root $\veps_1-\veps_2$. 
\end{m-remark}


\subsection{When is $G/P$ a $1$-splitting variety?}\label{ssct:1split}

To settle this question, we need some notations. 
\begin{longtable}{>{\raggedleft}p{5em} >{\raggedright}p{0.75\linewidth}}
	$\cal X^*(T):=$
& 
	the group of characters of $T$; similarly for $B,P,G$, etc.;
\tabularnewline 
	$(\,\cdot\,,\cdot\,)$
& 
	the $\Weyl(G)$-invariant scalar product on $\cal X^*(T)_{\mbb Q}$;
\tabularnewline 
	$\Psi:=$ 
& 
	the roots of $G$,\qquad $\Delta\subset\Psi$ the simple roots; 
\tabularnewline 
	$\L:=$
&
	the weights 
	$\{\omega\in\cal X^*(T)_{\mbb Q}\mid (\omega,\beta^\vee)\in\mbb Z,\;\forall\,\beta\in\Delta\}$; $\beta^\vee:=\frac{2\beta}{(\beta,\beta)}$; 
\tabularnewline 
&
	$\{\omega_\alpha\}_{\alpha\in\Delta}$ the fundamental weights, 
	that is $(\omega_\alpha,\beta^\vee)={\rm Kronecker}_{\alpha\beta}$;
\tabularnewline 
	$\L_+:=$
&
	the dominant weights 
	$\{\omega\in\L\mid (\omega,\beta)\ges0,\;\forall\,\beta\in\Delta\}$; 
\tabularnewline 
	$\lran{I}:=$ 
& 
	the vector space generated by $I\subset\Delta$;
\tabularnewline 
	$\L(I):=$
&
	$\L\cap\underset{\alpha\in I}{\bigcap}\alpha^\perp=\L\cap\lran{I}^\perp$, 
	the $I$-face of $\L$. 
\end{longtable}
\nit The parabolic subgroup $P$ corresponds to a subset $I\subset\Delta$ 
(cf. \cite[Section 8.4]{sp}); we denote it by $P_I$. 
Its Weyl group $W_I$ is generated by the reflections $\tau_\alpha$, $\alpha\in I$. 

\begin{m-proposition}\label{prop:homog-1split}
The homogeneous space $X=G/P_I$ is $1$-splitting (cf. \ref{def:1split}) if and only if 
there is no simple root perpendicular to $\lran{\alpha\mid\alpha\in I}$. 
Equivalently, define 
$$\tld I:=
I\cup\{
\beta\in\Delta\sm I\mid \beta\text{ is adjacent to some }\alpha\in I
\text{ in the Dynkin diagram of }G\}.$$
Then $G/P_I$ is $1$-splitting if and only if $\tld I=\Delta$. 
(For Dynkin diagrams, see \cite[\S 9.5]{sp}.)
\end{m-proposition}

\begin{proof}
We may assume that $G$ is simply connected, so line bundles on $X$ correspond 
to characters of $P_I$. For $\chi\in\cal X^*(P_I)$, define 
$\eL_\chi:=(G\times\mbb C)/P_I$, where $(g,z)\sim(gp^{-1},\chi(p)z)$. 
By the Borel-Weil-Bott theorem \cite{bott-homog,demazure}, 
$H^1(G/P_I,\eL_\chi)= H^1(G/B,\eL_\chi)\neq 0$ if and only if 
there is $\beta\in\Delta$ and $\chi_+\in\L_+$ such that 
$$\chi=\tau_\beta(\chi_++\rho)-\rho=\tau_\beta(\chi_++\beta)
\;\Leftrightarrow\;
\chi_+=\tau_\beta(\chi+\beta).$$
Let $T_\beta$ be the transformation $\chi\mt \tau_\beta(\chi+\beta)$; 
it is the reflection in the plane orthogonal to $\beta$, passing through $-\frac{\beta}{2}$.
The pull-back of $\eL_\chi$ to $G/B$ corresponds to the image of $\chi$ by 
$\varpi_I:\cal X^*(P_I)\to\cal X^*(B)=\cal X^*(T)$, 
so $G/P_I$ is $1$-splitting if and only if 
$$\L_+\cap T_\beta(\Img(\varpi_I))=\emptyset,\;\forall\,\beta\in\Delta.$$
But $\cal X^*(P_I)=\cal X^*(\Levi(P_I))$, so $\Img(\varpi_I)\subset\cal X^*(T)$ 
consists of the $W_I$-invariant elements. Since $W_I$ is generated by the reflections 
$\tau_\alpha$ in the hyperplanes $\alpha^\perp$, $\alpha\in I$, we deduce that 
$\Img(\varpi_I)=\L(I)$, so we must have 
$$\L_+\cap T_\beta\big(\,\L(I)\,\big)=\emptyset,\;\forall\,\beta\in\Delta.$$
For $\beta\in I$, this condition is automatically satisfied: 
$\L(I)\subset\beta^\perp$, so 
$$T_\beta\big(\,\L(I)\,\big)\subset\{(\beta,\cdot)<0\}
\text{ and  }
\L_+\subset\{(\beta,\cdot)\ges0\}.$$ 
For $\beta\in\Delta\sm I$, let $\beta_{\perp}$ be the component of $\beta$ 
on $\lran{I}$ with respect to the orthogonal decomposition 
$\lran{\cal X^*(T)}=\lran{I}\oplus\lran{I}^\perp$.  Then $\beta_{\perp}$ is also the 
orthogonal projection of $0$ to the affine space $\beta+\lran{I}^\perp$; hence $\L_+$ 
and $T_\beta\big(\,\L(I)\,\big)=\tau_\beta\big(\beta+\L(I)\big)$ 
are disjoint if and only if they are on different sides of the hyperplane 
$\lran{\tau_\beta\beta_\perp}^\perp$: 
\begin{equation}\label{eq:gp-split}
{\rm(i)}\;(\beta_\perp,\beta)>0,
\quad{\rm(ii)}\;(\tau_\beta\beta_\perp,\omega_\alpha)
=(\beta_\perp,\tau_\beta\omega_\alpha)\les0,
\;\forall\,\beta\in\Delta\sm I,\;\forall\,\alpha\in\Delta.
\end{equation} 
\nit\textit{Claim} The inequality (ii) is automatically satisfied. 

\begin{longtable}{rl}
\hspace{-1\parindent}\textit{Case $\alpha\neq\beta$}: 
&
$\tau_\beta\omega_\alpha=\omega_\alpha$.
\tabularnewline 
& 
-- Assume $\alpha\not\in I$. Then holds 
$\omega_\alpha\in\L_+(I)\;\Rightarrow\;(\beta_\perp,\omega_\alpha)=0$. 
\tabularnewline 
&
\begin{minipage}[t]{0.8\linewidth}
-- Assume $\alpha\in I$. 
Since any two vectors in $I$ make an angle of at least $90^\circ$ (they are simple roots), 
and $(-\beta_\perp,c)=-(\beta,c)\ges 0$, for all $c\in I$, we deduce that $-\beta_\perp$ 
is in the cone $\underset{c\in I}{\sum}\,\bbR_{\ges 0}c$. Thus holds:
\\[.5ex] \centerline{
$
-\beta_{\perp}=k_\alpha\alpha+\underset{c\in I\sm\{\alpha\}}{\sum}k_cc,\;k_c\ges 0
\quad\Rightarrow\quad
(\beta_{\perp},\omega_\alpha)=-k_\alpha(\alpha,\omega_\alpha)\les 0. 
$
}\smallskip
\end{minipage}
\tabularnewline
\textit{Case $\alpha=\beta$}: 
&
\begin{minipage}[t]{0.8\linewidth}
$\omega_\beta\in\L_+(I) 
\;\Rightarrow\; 
\beta_\perp\perp\omega_\beta
\;\Rightarrow\; 
(\beta_\perp,\tau_\beta\omega_\beta)=
(\beta_\perp,\omega_\beta-\beta)=
-(\beta_\perp,\beta)\les 0.$ 
For the last step: the angle between a vector and its projection to any plane 
is at most $90^\circ$. \vspace{1ex}
\end{minipage}
\end{longtable}

\nit Hence the only relevant condition in \eqref{eq:gp-split} is the first one. 
However, $(\beta_\perp,\beta)\ges 0$ from the very construction, so we must 
eliminate the case $(\beta_\perp,\beta)= 0$. This happens precisely when 
$0\in\beta+\lran{\L(I)}
\;\Leftrightarrow\;
\beta\in\L(I)
\;\Leftrightarrow\;
\beta\perp\alpha,\;\forall\,\alpha\in I.$
\end{proof}

\begin{m-corollary}\label{cor:iterate}
Assume $T\subset B^-\subset P_I\subset G$, 
the variety $X=G/P_I$ is $1$-splitting, and 
\\[1ex] \centerline{
$2\,{\cdot}\#I\ges1+\#\Delta$\quad (here $\#$ stands for the cardinality).
}\\[1ex] 
Then there is $\l:G_m\to T$ such that the source $Y$ of the action has the properties: 
$$\begin{array}{l}
\phantom{t}{\rm(i)}\;
Y=G(\l)/G(\l)\cap P_I;
\quad
{\rm(ii)}\;
Y\text{ is $1$-splitting};
\\[1ex] 
{\rm(iii)}\;
G(\l)\cap P_I\text{ corresponds to the simple roots }I\sm\{\alpha_0\},\;\alpha_0\in I.
\end{array}$$
\end{m-corollary}

\begin{proof}
Denote by $\Psi_I$ the roots generated by $I$; the roots of 
the unipotent radical $R_u(P_I)\subset P_I$ are $\Psi^-\sm\Psi^-_I$ 
(cf. \cite[Theorem 8.4.3]{sp}). 
For an arbitrary $\alpha_0\in I$, we consider a 1-PS $\l$ such that 
$\lran{\alpha_0,\l}>0,$ $\,\lran{c,\l}=0,\;\forall\,c\in\Delta\sm\{\alpha_0\}.$ 
One can easily check the following:
$$
\begin{array}{rl}
\mbox{\rm negative roots of }G(\l)&=\big\{
\underset{c\in\Delta\sm\{\alpha_0\}}{\sum}\kern-1ex
k_cc\in\Psi\mid k_c\les 0\big\} \supset\Psi_I^-, 
\\[2.5ex] 
\Lie(G)_\l^-&=\big\{
k_0\alpha_0+\underset{c\in\Delta\sm\{\alpha_0\}}{\sum}\kern-1exk_cc
\mid k_0<0,\;k_c\les0\;\mbox{\rm for }c\neq\alpha_0
\big\}\subset\Psi^-.
\end{array}
$$
It follows that $\Lie(G)_\l^-\subset\Psi^-\sm\Psi^-_I$, which are the roots of 
$R_u(P_I)$, which implies that $\Lie(G)_\l^-\subset\Lie(P_I)$, as desired. 
It remains to prove that $\alpha_0$ can be chosen in such a way that 
$Y$ is $1$-splitting, that is:
$$
\not\kern-2pt\exists\,\beta\in \Delta\sm I\text{ such that }
\beta\text{ is adjacent only to }\alpha_0.
$$
If such an $\alpha_0\in I$ does not exist, then for all $\alpha\in I$ there is 
$\beta\in\Delta\sm I$ adjacent only to $\alpha$; this yields an injective function 
$I\to\Delta\sm I$, so $\#I\les\#(\Delta\sm I)$, which contradicts the hypothesis.
\end{proof}


\section{Application: splitting criteria for vector bundles on Grassmannians} \label{sct:grass}

The Grassmannian plays a central role because it is a homogeneous variety, and is the 
`universal target' for pairs $(X,\eN)$ consisting of a variety and a globally generated 
vector bundle on it. 
Hence the results of both the sections \ref{sct:split-0loci}, \ref{sct:split-homog} apply. 
In this section we obtain splitting criteria for vector bundles on the isotropic 
(symplectic and orthogonal) Grassmannians. 
The degenerate case, when the bilinear form has kernel, is also included, to demonstrate 
that theorem \ref{thm:split} is not restricted only to the situations discussed in sections 
\ref{sct:0vb}, \ref{sct:fixed}. All the computations involve two stages: first we compute 
the positivity of various subvarieties; second, we deduce splitting criteria by restricting 
vector bundles on the ambient space to them. 

Cohomological splitting criteria have been obtained in \cite{ott,ma,ar-ma,ma-oe}; 
however, they involve a large number of conditions. 
The results below are interesting for their simplicity: indeed, 
the problem of deciding the splitting of a vector bundle on a Grassmann variety, 
which is a high dimensional object, is reduced to the splitting along a (very) 
low dimensional subvariety. 
Throughout this section, $W$ stands for a $w+1=\nu+u+1$-dimensional vector space.


\subsection{The Grassmannian of linear subspaces}\label{ssct:gl}

This case is discussed in \cite{hal}, where is proved that things are as good as possible, 
without any genericity assumptions. 
\begin{m-theorem}\label{thm:split-gl}
The vector bundle $\crl V$ on $\Grs(u;\bbC^{w})$, $u\ges 2,\,w\ges u+2$, splits 
if and only if its restriction to \emph{an arbitrary} $\Grs(2;\bk^4)\subset\Grs(u;\bk^{w})$ 
does so.
\end{m-theorem}
This is in perfect analogy with Horrocks' criterion. However, the proof uses 
a (fortunate) cohomology vanishing, and can not be extended directly. 
As a warm-up, let us see compute the positivity of a `smaller' Grassmannian 
in `larger' one. 
\begin{m-example}\label{expl:Grs}
Let $X=\Grs(u+1;\mbb C^{w+1})$ be the Grassmannian of $(u+1)$-dimensional 
subspaces of $\mbb C^{w+1}$. Let $e_0,e_1,\ldots$ be the standard basis of 
$\mbb C^{w+1}$. The $1$-PS 
$$
\l:G_m\to\GL(w+1),\quad\l(t):=\diag(t^{-1},1,\ldots,1)
$$ 
acts on $X$ with source 
$Y=\{U\in X\mid e_0\in U\}=\Grs(u;\mbb C^{w+1}/\lran{e_0})$. 
For $U\in Y$, the normal bundle is 
$\eN_{Y/X,U}=\Hom(\lran{s},\mbb C^{w+1}/U)$; 
since $e_0\in U$, $\l$ acts trivially on $\mbb C^{w+1}/U$, so it acts by scalar 
multiplication on the normal bundle. One can easily see that 
$$
X\sm Y^+=\{U\in X\mid U\subset\lran{e_1,\ldots,e_w}\}\cong\Grs(u+1;w),
$$ 
so $\Grs(u;\mbb C^{w})\subset\Grs(u+1;\mbb C^{w+1})\text{ is }u^\pos$ 
(cf. proposition \ref{prop:y-sink}). 
\end{m-example}


\subsection{The symplectic-isotropic Grassmannian}\label{ssct:sp}

Let $\omega$ be a skew-symmetric bilinear form on $W$, such that: 

-- if $\dim W$ is even, 
$\omega$ is non-degenerate (so $\omega$ is a usual symplectic form); 

-- if $\dim W$ is odd, $\dim\Ker(\omega)=1$ 
($\omega$ is symplectic on $W/\Ker(\omega)$ of dimension $w$). 

\nit Let $X:=\spGrs(u+1;W)$ be the variety of $\omega$-isotropic, 
$(u+1)$-dimensional subspaces of $W$. It is a Fano variety with 
$$
\dim(X)=\frac{(u+1)(2w-3u)}{2}.
$$ 
If $\vphi:\spGrs(u+1;W)\to\Grs(u+1;W)$ stands for the natural embedding, then 
$\eO_{X}(1)$ is the pull-back of the corresponding line bundle on the Grassmannian. 

Denote by $G:=\Sp_{[(w+1)/2]}$ the symplectic group ($[\,\cdot\,]$ stands for 
the integral part). If $\dim W$ is even, $X$ is homogeneous for the $G$-action: 
$X=G/P$, where $P$ is the stabilizer of the flag 
$$
\big\{
(1,\ldots,\underset{u+1}{0},\ldots,0\mid 0,\ldots,0),\ldots,
(0,\ldots,\underset{u+1}{1},0,\ldots,0\mid0,\ldots,0)
\big\}
\subset\mbb C^{\frac{w+1}{2}}\oplus\mbb C^{\frac{w+1}{2}}.
$$
If $\dim W$ is odd, $X$ has two $G$-orbits: the open orbit of subspaces 
which intersect $\Ker(\omega)$ trivially, and the closed orbit of subspaces containing 
$\Ker(\omega)$. 
\begin{m-lemma}\label{lm:sp-pic}
If $w\ges 2u+1+\dim\Ker(\omega)$, then $\Pic(X)=\mbb Z\cdot\eO_X(1).$
\end{m-lemma}

\begin{proof}
For $\omega$ non-degenerate, this is clear. For $\Ker(\omega)=\lran{s_0}$, 
$$
\spGrs(u+1;W)=\{U\mid s_0\in W\}\cup\{U\mid s_0\not\in U\}.
$$
The first term is a subvariety, isomorphic to $\spGrs(u;W/\lran{s_0})$, of codimension 
$w-2u\ges 2$; thus $\Pic(X)$ is isomorphic to the Picard group of the open stratum. 
The morphism 
$$
\{U\mid s_0\not\in U\}\to\spGrs(u+1,W/\lran{s_0}), \quad 
[U\subset W]\mt\big[\lran{s_0}+U/\lran{s_0}\subset W/\lran{s_0}\big]
$$
is an affine space fibration, and the base has cyclic Picard group.
\end{proof}

The quotient bundle $\eN$ and the tautological bundle 
$\eU:=\Ker(W\otimes\eO_X\to\eN)$ on $X$ are the pull-back by $\vphi$ 
of their counterparts on the Grassmannian. 
An element $s\in W\sm\Ker(\omega)$ defines a section in $\eN$, whose zero set is 
the `smaller' isotropic Grassmannian: 
\begin{equation}\label{eq:Ws}
\begin{array}{r}
\{U\in\spGrs(u;W)\mid s\in U\}=
\{U\in\spGrs(u;W)\mid s\in U\subset U^\perp\subset \lran{s}^\perp\}
\\[1ex] 
=\spGrs(u;\lran{s}^\perp/\lran{s})=\spGrs(u+1;W)\cap \Grs(u;W/\lran{s}),
\\[1.5ex]
\text{where}\;\lran{s}^\perp:=\{t\in W\mid\omega(s,t)=0\}.\hfill\null
\end{array}
\end{equation}
An element $\si\in W^\vee\sm\{0\}$ determines a section in $\eU^\vee$, with zero locus  
\begin{equation}\label{eq:Wsi}
\spGrs(u+1,\si^\perp)=\{U\in\spGrs(u+1;W)\mid U\subset\si^\perp\},
\;\text{where}\;
\si^\perp:=\Ker(\si).
\end{equation}
In particular, $s\in\big(W/\Ker(\omega)\big)\sm\{0\}$ determines 
$\si_s(\cdot):=\omega(s,\cdot)\in W^\vee$; in this case, $\si_s^\perp=\lran{s}^\perp$. 
For $\omega$ non-degenerate, $s\mt\si_s$ defines an isomorphism $W\to W^\vee$; 
however, if $\dim\Ker(\omega)=1$, the image of this map is a hyperplane in $W^\vee$. 
In the latter case, for generic $\si\in W^\vee$, $\omega|_{\si^\perp}$ is 
non-degenerate, so $\si^\perp$ is a symplectic subspace of $W$. In general, it holds 
\begin{equation}\label{eq:ker}
\dim\Ker(\omega|_{\si^\perp})=1-\dim\Ker(\omega), 
\;\text{for generic}\;\si.
\end{equation}

We start by explicitly computing the positivity of some subvarieties of $X$. Deliberately, 
we consider both zero loci of sections and sources of $G_m$-actions, to illustrate 
the general theory developed in the previous sections. 

\begin{m-lemma}\label{lm:sp-p-pos}
\nit{\rm(i)} 
For $s\in W^\vee\sm\{0\}$ and $u\ges 1$, 
$$
\spGrs(u+1;s^\perp)\subset\spGrs(u+1;W)
\text{ is }(w-2u-1)^\pos.\quad(\text{Recall that $s^\perp:=\Ker(s)$.})
$$

\nit{\rm(ii)} 
Assume that either $\Ker(\omega)=\lran{s}\subset W$, 
or $\omega$ is non-degenerate and $s\neq0$. Then 
$$
\spGrs(u;\lran{s}^\perp/\lran{s})\subset\spGrs(u+1;W)
\text{ is }{u}^\pos,\text{ for }u\ges 1.
$$

\nit{\rm(iii)} Decompose $W=W'\oplus W''$ into the sum of Lagrangian subspaces 
of dimension $(w+1)/2$. Then 
$
\Grs(u+1;W')\subset\spGrs(u+1,W)\text{ is }\big(\frac{w-2u-1}{2}\big)^\pos.
$ 
\end{m-lemma}

\begin{proof}
(i) In the diagram 
\begin{equation}\label{eq:spgrs1}
\xymatrix@R=1.5em@C=5em{
\spGrs(u+1;W)\ar@{-->}[r]_-{U\mt U\cap s^\perp}^-b\ar[d]^\vphi
&
\spGrs(u;s^\perp)\ar[d]^\vphi
\\ 
\Grs(u+1;W)\ar@{-->}[r]_-{U\mt U\cap s^\perp}^-{g}
&
\Grs(u;s^\perp),
}
\end{equation}
$g$ is undefined on $\Grs(u+1;s^\perp)$, and $b$ is undefined on 
$\{U\mid U\subset s^\perp\}=\spGrs(u+1;s^\perp)$. 
The blow-up of $\Grs(u+1;W)$ along $\Grs(u+1;s^\perp)$ resolves $g$, hence $b$. 
It remains to apply proposition \ref{prop:p}: the fibres of $b$ are at least 
$\frac{(u+1)(2w-3u)}{2}-\frac{u(2w-3u+1)}{2}=w-2u$ dimensional. 

\nit(ii) \textit{Case $\Ker(\omega)=\lran{s}$.}\quad 
In the diagram 
\begin{equation}\label{eq:spgrs2}\xymatrix@C=5em@R=1.5em{
\spGrs(u+1;W)\ar@{-->}[r]_-{U\mt (\lran{s}+U)/\lran{s}}^-b\ar[d]^\vphi
&
\spGrs(u+1;W/\lran{s})\ar[d]^\vphi
&
\kern-4em
\big(
\text{$\frac{\lran{s}+U}{\lran{s}}\subset\frac{W}{\lran{s}}$ is isotropic.}
\big)
\\ 
\Grs(u+1;W)\ar@{-->}[r]_-{U\mt (\lran{s}+U)/\lran{s}}^-{g}&\Grs(u+1;W/\lran{s}).
&
}
\end{equation}
$b$ is not defined on $\Grs(u;W/\lran{s})\cap\spGrs(u+1;W)=\spGrs(u;W/\lran{s})$. 
The blow-up of $\Grs(u;W/\lran{s})$ resolves $g$, hence $b$.  
Now apply \ref{prop:p} again: the general fibre of $b$ is $(u+1)$-dimensional. 
\smallskip 

\nit\textit{Case $W$ is symplectic.}\quad 
To settle this, we use the proposition \ref{prop:trans} in conjunction with 
\ref{prop:cdX-Y}. Decompose $W$ into a direct sum of Lagrangian subspaces  
$W=\bk^{(w+1)/2}_{\rm left}\oplus\bk^{(w+1)/2}_{\rm right},$ 
and assume that $s=(1,\ldots,0\,|\,0,\ldots,0)$. 
Then $Y:=\spGrs(u,\lran{s}^\perp/\lran{s})$ is the source of the $G_m$-action, 
corresponding to the $1$-PS: 
\begin{equation}\label{eq:l1}
\begin{array}{r}
\l:G_m\to\Sp_{(w+1)/2}(\bk),\quad
\l(t)=\diag\Big[t^{-1},\bone_{(w-1)/2},t,\bone_{(w-1)/2}\Big]
\end{array}
\end{equation}
The complement of the open BB-cell is 
$X\sm Y^+=\big\{U\in X\mid s\not\in\lim_{t\to 0}\l(t)U\big\}.$ 
Consider a basis in $U$ such that the corresponding column matrix is lower triangular; 
then one sees that 
$X\sm Y^+{=}
\,\{U\mid U\subset\bk^{(w-1)/2}_{\rm left}\oplus\bk^{(w+1)/2}_{\rm right}\}$, 
which is $\frac{(u+1)(2(w-1)-3u)}{2}$-dimensional, so 
$$
\cd(X\sm Y)=\dim X-(u+1).
$$ 
It remain to compute the amplitude of the normal bundle $\eN_{Y/X}$, which is 
isomorphic to the universal quotient bundle $\eN$ on $Y$. This is globally generated, 
so we can use the criterion \ref{prop:q2} to compute its ampleness. 
The homomorphism $\frac{W}{\lran{s}}\otimes\eO_Y\to\eN$ induces 
the surjective map $\mbb P(\eN^\vee)\to\mbb P^{w-1}$. By homogeneity, 
its fibres are isomorphic, hence $\eN_{Y/X}$ is $(\dim Y-u)$-ample. 
 
\nit(iii) 
We are going to apply the proposition \ref{prop:y-sink}. 
For a Lagrangian direct sum decomposition 
$W=W'\oplus W''=\bk^{(w+1)/2}_{\rm left}\oplus\bk^{(w+1)/2}_{\rm right},$ 
we observe that $Y:=\Grs(u+1,\bk^{(w+1)/2}_{\rm left})$ is the source of 
the $G_m$-action 
\begin{equation}\label{eq:l2}
\l:G_m\to\Sp_{(w+1)/2}(\bk),\quad
\l(t)=\diag\big[t^{-1}\bone_{(w+1)/2},t\bone_{(w+1)/2}\big].
\end{equation}
An easy computation yields: 
\\[.5ex]\centerline{
$\begin{array}{rl}
T_{X,U}
&=
\{h\in\Hom\big(U,W/U\big)\mid \omega(u',hu'')+\omega(hu',u'')=0,\forall\,u',u''\in U\}
\\[.5ex] 
&\cong\Hom(U,U^\perp/U)\oplus\Hom^{\rm symm}(U,U^\vee), 
\\[1ex] 
T_{Y,U}
&
=\Hom(U,W'/U), 
\\[1ex] 
\eN_{Y/X, U}
&=\Hom(U,U^\perp/W')\oplus\Hom^{\rm symm}(U,U^\vee). 
\end{array}$
}\\[.5ex]
On both summands $\l$ acts diagonally, with weight $t^2$. Furthermore, 
one can see that 
$$
X\sm Y^+=\{U\mid \Ker(\pr:U\to\bk^{(w+1)/2}_{\rm left})\neq 0\}.
$$ 
The minimal degeneration is when $\Ker(\pr)$ is one dimensional; the corresponding 
stratum maps onto $\mbb P(\bk^{(w+1)/2}_{\rm right})$, with fibres isomorphic to 
$\spGrs(u;\bk^{w-1})$. It follows that $\dim X-\dim(X\sm Y^+)=\frac{w-2u+1}{2}$. 
\end{proof}

\begin{m-remark}\label{rmk:sp-bad}
\nit{\rm(i)} 
Sommese's criterion \ref{prop:q2} implies that the quotient bundle on $\spGrs(u+1;W)$ 
is $q$-ample, for $q=\dim\spGrs(u+1;W)-(u+1)$. 
Hence $\Grs(u;\lran{s}^\perp/\lran{s})$ is $p^\pos$, with 
$$
p=(u+1)-\codim\,\spGrs(u;\lran{s}^\perp/\lran{s})=2u-w+1<0.
$$
This shows that this test is week compared with proposition \ref{prop:p}. 

\nit{\rm(ii)} The conclusion of \ref{lm:sp-p-pos}(i), for $\omega$ non-degenerate, 
can not be obtained by using a $1$-PS of $G$, because $\spGrs(u+1;s^\perp)$ is not 
homogeneous. 
In contrast, it is not clear how to prove \ref{lm:sp-p-pos}(iii) by using 
the proposition \ref{prop:p}. 

\nit{\rm(iii)} At \ref{lm:sp-p-pos}(ii), for $\omega$ degenerate, the section $s$ 
(with $\lran{s}=\Ker(\omega)$) is \emph{neither generic nor regular} 
(transverse to zero); thus we really use the generality in the subsection \ref{ssct:x-b}. 

For $\omega$ non-degenerate, the computation \ref{lm:sp-p-pos}(ii) is one 
of the most technical ones. None of the other methods 
(propositions \ref{prop:p} and \ref{prop:y-sink}) yield the $u^\pos$ positivity. 
Proposition \ref{prop:trans} readily implies that 
$\spGrs(u;w-1)\subset\spGrs(u+1;w+1)$ is ${\min(u,w-2u-1)}^\apos$, 
but this is too weak to obtain splitting criterion for $\spGrs(u;2u)$. 
\end{m-remark}

\begin{m-proposition}\label{prop:split-sp}
Let $\omega$ be a skew-symmetric bilinear form on $W$, $\dim W=w+1$, and 
$\kappa:=\dim\Ker(\omega)\les 1$. Let $X:=\spGrs(u+1;W)$, and $\crl V$ be 
an arbitrary vector bundle on $X$. 
In the cases enumerated below, $\crl V$ splits if and only if $\crl V_{Y_s}$ splits. 

\nit{\rm(i)} $Y_s:=\spGrs(u+1;s^\perp)$, with $s\in W^\vee$ very general, and  
\\[1ex] \centerline{
$w+1>2(u+1)+\kappa$, so $w\ges 2u+3+\kappa$, and $u\ges1$.
}\\[1ex] 
\nit{\rm(ii)} $Y_s:=\spGrs(u;\lran{s}^\perp\!/\lran{s})$, 
with $s\in W{\sm}\,\Ker(\omega)$ very general, and $w\ges 2u+1+\kappa$, $u\ges 2$.
\end{m-proposition}

\begin{proof}
In both cases, we verify \ysx and apply theorem \ref{thm:split} directly. 

\nit(i)  Take $S\subset W^\vee\sm\{0\}$ open subset such that \eqref{eq:ker} holds, 
and $\cY\subset S\times X$ be the zero locus of the universal section 
in ${(\pr^{S\times X}_X)}^*\eU^\vee$. 
The morphisms $\cY\srel{\pi}{\to}S$ and $\cY\srel{\rho}{\to} X$ 
are both open, and $\pi^{-1}(s)=Y_s$ is projective, connected; also, 
$\dim\Ker(\omega|_{s^\perp})=1-\kappa$.  

\nit-- We claim that, for all $o,s,t\in S$,  
$$
\begin{array}{l}
Y_{os}=
\{U\in X\mid U\subset\lran{o,s}^\perp\}\cong\spGrs(u+1;\lran{o,s}^\perp),
\text{ and}
\\[1ex] 
Y_{ost}=\spGrs(u+1;\lran{o,s,t}^\perp)\cong\spGrs(u+1;\lran{o,s,t}^\perp)
\end{array}
$$
are connected and non-empty. 
The connectedness is clear, since they are quasi-homogeneous, with finitely many orbits, 
for actions of appropriate subgroups of $\Sp(W)$. 
We verify that $Y_{ost}\neq\emptyset$ (thus \arm and \notr are satisfied), that is 
$$
\exists\,U\in X\text{ such that }U\subset\lran{o,s,t}^\perp.
$$
Let $\omega_{ost}$ be the restriction of $\omega$ to $\lran{o,s,t}^\perp$; 
then $\dim\Ker(\omega_{ost})=1-\kappa$. We verify: 
$$
\dim\lran{o,s,t}^\perp/\,\Ker(\omega_{ost})\ges 
2\big((u+1)-\dim\Ker(\omega_{ost})\big)
\;\Leftrightarrow\;
w+\dim\Ker(\omega_{ost})\ges 2u+4.
$$
The last inequality is fulfilled, by hypothesis. 

\nit-- The diagram \ysx\hskip-1ex(Pic) consists of isomorphisms 
(cf. lemma \ref{lm:sp-pic}).

\nit-- Finally, $Y_s\subset X$ is $(w-2u-1)^\pos$ 
(cf. lemma \ref{lm:sp-p-pos}(i)), so at least $2^\pos$. 
\smallskip 

\nit(ii) Let $S:=W\sm\Ker(\omega)$, 
and $\cY$ be the zero locus of the universal section in ${(\pr^{S\times X}_X)}^*\eN$. 
The situation is analogous, with a few differences. Indeed, for $o,s\in S$ holds: 
$$
Y_{os}\neq\emptyset\;\Leftrightarrow\;
o\perp s, 
\text{ so generically }Y_{os}=\emptyset.
$$
For $o\in S$, we check the properties \ysx for 
$$
S(o)=\{s\in S\mid Y_{os}\neq\emptyset\}=S\cap\lran{o}^\perp.
$$

\nit -- The diagram (Pic) consists of isomorphisms (cf. lemma \ref{lm:sp-pic}). 

\nit-- The condition \arm\kern-3pt, that is $\rho(\cY_{S(o)})=X$, is the following: 
$$
\forall\,U\in X\;\exists\,s\in S\cap\lran{o}^\perp\;\exists\,V\in X\text{ such that }
s\in U,\;\lran{o,s}\subset V.
$$
Indeed, $\dim(U\cap\lran{o}^\perp)\ges u$, so $U\cap\lran{o}^\perp\neq 0$. 
Take $s\in(U\cap\lran{o}^\perp)\sm\Ker(\omega)$, non-zero; 
then $\lran{o,s}\subset W$ is an isotropic subspace. 
There exists $V\in X$ containing it, as $u\ges 2$. 

\nit-- The condition \notr reads:
$$
[s\perp o,\,t\perp o,\,s\perp t]
\;\Rightarrow\;\exists\,V\in X,\;\lran{o,s,t}\subset V.
$$
The left hand side implies that $\lran{o,s,t}\subset W$ is isotropic subspace; 
then $V$ exists, since $u\ges2$. 

\nit-- Finally, by lemma \ref{lm:sp-p-pos}(ii), $Y_s\subset X$ 
is $u^\pos$, so at least $2^\apos$. 
\end{proof}

\begin{m-theorem}\label{thm:split-sp}
Let $\omega$ be a skew-symmetric bilinear form on $\bk^{w}$, with 
$\kappa=\dim\Ker(\omega)\les1$. We consider the isotropic Grassmannian 
$X=\spGrs(u;\bk^{w})$, with $u\ges 2$, and an arbitrary vector bundle $\crl V$ 
on it. Then $\crl V$ splits if and only if it does so along a very general subvariety 
$Y\cong\spGrs(2;4+\kappa)$ of $X$. 
\end{m-theorem}

Note that $\spGrs(2;4)$, the Lagrangian $2$-planes in $\bk^4$, 
is isomorphic through the Pl\"ucker embedding to the $3$-dimensional quadric. 

\begin{proof} 
By applying repeatedly the first part of previous proposition 
(after replacing $w+1\rightsquigarrow w$ and $u+1\rightsquigarrow u$), we deduce  
that $\crl V$ splits if and only if $\crl V_Z$ splits, for some very general subvariety 
$Z\cong\spGrs(u,2u+\kappa)$ of $X$. Now apply the second part to deduce reduce 
the splitting problem from $Z$ to $Y\cong\spGrs(2,4+\kappa)$. 
(For this latter, the process can not be iterated anymore.) 
\end{proof}


\subsection{The orthogonal-isotropic Grassmannian}\label{ssct:o}

The situation is similar to the previous case: let $\beta$ be a symmetric, non-degenerate, 
bilinear symmetric form on $W$, and consider $X:=\oGrs(u+1;W)$ be the variety of 
isotropic $(u+1)$-dimensional subspaces of $W$; assume $u\ges 1$, $\dim W\ges 5$. 
(If $w+1=2(u+1)$, the full space of Lagrangian planes in $\bk^{2(u+1)}$ 
has two connected components, and we consider only one of them.) 
It is a homogeneous variety for $G=\SO(\beta)$, with 
$$
\dim X=\frac{(u+1)(2w-3u-2)}{2},\text{ and }\Pic(X)=\mbb Z\cdot\eO_X(1). 
$$ 
Similar arguments as before yield the following.

\begin{m-lemma}\label{lm:o-p-pos}
\nit{\rm(i)} If $u\ges 1$, $w\ges 2u+2$, then 
$$
\oGrs(u+1;\lran{s}^\perp)\subset\oGrs(u+1;W)\text{ is }(w-2u-2)^\pos. 
$$
The form $\beta|_{\lran{s}^\perp}$ is non-degenerate if and only if 
$s$ is a non-isotropic vector.

\nit{\rm(ii)} 
Let $s\in W$ be isotropic, $u\ges 2$. Then 
$\oGrs(u;\lran{s}^\perp/\lran{s})\subset\oGrs(u+1;W)$ is 
\\[1ex] \centerline{
$\biggl\{\begin{array}{cl}
(u-1)^\pos&\text{ for }w=2u+1,
\\[1ex] 
u^\pos&\text{ for }w\ges2u+2.
\end{array}\biggr.$
}
\end{m-lemma}

\begin{proof} 
(i) 
We repeat the argument in the lemma \ref{lm:sp-p-pos}(i), 
and use that $\dim\oGrs(u+1;W)-\dim\oGrs(u;\lran{s}^\perp)=w-2u-1$. 

\nit(ii) We are going to apply the proposition \ref{prop:y-sink}. 
Decompose $W=\mbb C^{(w+1)/2}\oplus\mbb C^{(w+1)/2}$ into the sum 
of two Lagrangian subspaces and consider the 1-PS: 
$$
\l:G_m\to\SO_{(w+1)/2},\quad
\l(t)=\big[t^{-1},\bone_{(w-1)/2},t,\bone_{(w+1)/2}\big].
$$ 
The source of $\l$ is $Y=\{U\mid s:=(1,0,\ldots,0)\in U\}$, 
and a simple computation shows that, for all $U\in Y$, holds: 
\\[.5ex]\centerline{
$\begin{array}{rl}
T_{X,U}
&=
\{h\in\Hom\big(U,W/U\big)\mid \beta(u',hu'')+\beta(hu',u'')=0,\forall\,u',u''\in U\}
\\[.5ex] 
&\cong\Hom(U,U^\perp/U)\oplus\Hom^{\rm anti-symm}(U,U^\vee),
\qquad\text{(Note that $h(s)\in\lran{s}^\perp$.)}
\\[1ex] 
T_{Y,U}
&
\cong\Hom(U/\lran{s},U^\perp/U)\oplus
\Hom^{\rm anti-symm}\big(U/\lran{s},(U/\lran{s})^\vee\big), 
\\[1ex] 
\eN_{Y/X, U}
&=\Hom(\lran{s},\lran{s}^\perp/U). 
\end{array}$
}\\[.5ex]
It follows that $\l$ acts by scalar multiplication on $\eN_{Y/X}$, with weight $t$. 
Furthermore, the complement of the open BB-cell is 
$X\sm Y^+=\big\{U\in X\mid s\not\in\underset{t\to 0}{\lim}\l(t)U\big\}.$ 
By taking a basis in $U$ such that the corresponding column matrix 
is lower triangular, one can see that  
$$
X\sm Y^+=\{U\mid U\subset W':=\mbb C^{(w-1)/2}\oplus\mbb C^{(w+1)/2}\}.
$$ 
Note that $\beta|_{W'}$ has a $1$-dimensional kernel; denote it by $\lran{s'}$. 
If $w\in\{2u+1,2u+2\}$, one can see that $s'\in U$ for all $U\in X\sm Y^+$. 
By using this remark one finds that $\dim X-\dim(X\sm Y^+)$ equals: 
$u$ for $w=2u+1$, and $u+1$ for $w\ges 2u+2$. 
\end{proof}

\begin{m-proposition}\label{prop:split-o}
Let $\beta$ be a non-degenerate symmetric bilinear form on $W$, with $\dim W=w+1$. 
Let $X:=\oGrs(u+1;W)$, and $\crl V$ be an arbitrary vector bundle on $X$. 
In the cases enumerated below, $\crl V$ splits if and only if $\crl V_{Y_s}$ splits. 

\nit{\rm(i)} $Y_s:=\oGrs(u+1;\lran{s}^\perp)$, 
with $s\in W$ very general, non-isotropic, and $w\ges 2u+4$. 

\nit{\rm(ii)} $Y_s:=\oGrs(u;\lran{s}^\perp/\lran{s})$, 
with $s$ very general isotropic, and $w\ges 2u+1$, $u\ges 3$. 
\end{m-proposition}

\begin{proof}
Again we check \ysx and apply \ref{thm:split} directly. 
Let $Q_\beta:=\{s\in W\mid\beta(s,s)=0\}$ be the isotropic cone.

\nit(i) 
Take $S:=W\sm Q_\beta$ and let $\cY\subset S\times X$ be the universal family, 
with $Y_s=\oGrs(u+1,\lran{s}^\perp)$. 

\nit-- For  $o,s,t\in S$, 
$$
\begin{array}{l}
Y_{os}=\{U\in X\mid U\subset\lran{o,s}^\perp\}, \quad 
Y_{ost}=\{U\in X\mid U\subset\lran{o,s,t}^\perp\}.
\end{array}
$$
They are quasi-homogeneous for the action of appropriate subgroups of $\SO(\beta)$, 
thus connected. As $w\ges 2u+4$, we deduce $\dim\lran{o,s,t}^\perp\ges2\dim U$, 
so $Y_{ost}$ is always non-empty; in particular, \arm and \notr are satisfied. 

\nit-- Finally, the diagram (Pic) consists of isomorphisms.
\smallskip

\nit(ii) Here we choose $S:=Q_\beta\sm\{0\}$; the situation is similar 
to \ref{prop:split-sp}(ii). For $o,s\in S$, 
$$
Y_{os}\neq\emptyset\;\Leftrightarrow\;s\in\lran{o}^\perp\cap S;\quad
\text{ let }S(o):=\lran{o}^\perp\cap S.
$$
We check the conditions \ysx for $\cY_{S(o)}$. 

\nit-- $W_{os}:=\lran{o,s}^\perp/\lran{o,s}$ has an induced non-degenerate 
orthogonal form, so $Y_{os}$ is connected. 

\nit-- The diagram (Pic) consists of isomorphisms.

\nit-- The condition \arm\kern-3pt, that is $\rho(\cY_{S(o)})=X$ is:
$$
\forall\,U\in X\;\exists\,s\in S(o)\;\exists\,V\in X\text{ such that }
s\in U,\;\lran{o,s}\subset V.
$$
Indeed, take $s\in U\cap\lran{o}^\perp$ non-zero, 
so $\lran{o,s}\subset W$ is isotropic; now take any $V$ containing it. 

\nit-- The condition \notr reads:
$\big[o\perp s,\;o\perp t,\;s\perp t\big]\;\Rightarrow\;Y_{ost}\neq\emptyset.$ 
\newline 
Indeed, $\lran{o,s,t}\subset W$ is isotropic, so there is $U\in X$ containing 
it because $u+1\ges 3$. 

\nit-- By lemma \ref{lm:o-p-pos}(ii), $Y_s\subset X$ is at least $2^\pos$. 
\end{proof}

\begin{m-theorem}\label{thm:split-o}
Let $\omega$ be a non-degenerate bilinear form on $\bk^{w}$. 
We consider the isotropic Grassmannian $X=\spGrs(u;\bk^{w})$, with $u\ges 3$,  
and an arbitrary vector bundle $\crl V$ on it. 
Then $\crl V$ splits if and only if it $\crl V_Y$ splits, with $Y$ very general, where:
\begin{itemize}
	\item $Y\cong\oGrs(3,\bk^6)$, if $w=2u$;
	\item $Y\cong\oGrs(3,\bk^7)$, if $w=2u+1$;
	\item $Y\cong\oGrs(3,\bk^8)$, if $w\ges 2u+2$. 
\end{itemize}
\end{m-theorem}

\begin{proof}
Assume that $w\ges 2u+2$. Then the first part of previous proposition implies 
(after replacing $w+1\rightsquigarrow w$ and $u+1\rightsquigarrow u$) 
that $\crl V$ splits if and only if $\crl V_Z$ splits, for some very general subvariety 
$Z\cong\spGrs(u,2u+2)$ of $X$. There remain three possibilities: 
$\oGrs(u,2u),\oGrs(u,2u+1),\oGrs(u,2u+2)$. 
The theorem follows now from the second part of the proposition.  
\end{proof}

The somewhat non-uniform formulation of the theorem, compared to \ref{thm:split-sp}, 
is due to the lack of sufficient positivity of 
$Y=\oGrs(u;2u+1)\subset\oGrs(u,2u+2)=X$, which is only $1^\pos$. 



\begin{thebibliography}{ooo}

\bibitem{ag} \textrm{A.~Andreotti, H.~Grauert}, 
\textit{Th\'eor\`eme de finitude pour la cohomologie des espaces complexes}. 
Bull. Soc. Math. France \textbf{90} (1962) 193--259.

\bibitem{ar} \textrm{D.~Arapura},
\textit{Partial regularity and amplitude}. 
Amer.~J.~Math. \textbf{128} (2006), 1025--1056. 

\bibitem{ar-ma} \textrm{E.~Arrondo, F.~Malaspina}, 
\textit{Cohomological characterization of vector bundles on Grassmannians of lines}. 
J.~Algebra \textbf{323} (2010), 1098--1106.

\bibitem{at} \textrm{M.~Atiyah}, 
\textit{On the Krull-Schmidt theorem with applications to sheaves}. 
Bull. Soc. Math. France \textbf{84} (1956), 307--317.

\bibitem{ba} \textrm{P.~Bakhtary}, 
\textit{Splitting criteria for vector bundles on higher dimensional varieties}. 
Pacific J.~Math. \textbf{252} (2011), 19--29.

\bibitem{bb} \textrm{A.~Bialynicki-Birula}, 
\textit{Some theorems on actions of algebraic groups}. 
Ann. Math. \textbf{98} (1973), 480--497.

\bibitem{bott-lefschetz} \textrm{R.~Bott}, 
\textit{On a theorem of Lefschetz}. 
Michigan Math.~J. \textbf{6} (1959), 211--216.

\bibitem{bott-homog} \textrm{R.~Bott}, 
\textit{Homogeneous vector bundles}. 
Ann. Math. \textbf{66} (1957), 203--248.

\bibitem{dps} \textrm{J.-P.~Demailly, T.~Peternell, M.~Schneider}, 
\textit{Holomorphic line bundles with partially vanishing cohomology}, 
in Proceedings of the Hirzebruch 65 Conference on Algebraic Geometry. 
Bar-Ilan Univ, Ramat Gan, 1996, 165--198.

\bibitem{demazure} \textrm{M.~Demazure},
\textit{A very simple proof of Bott's theorem}.
Invent. Math. \textbf{33} (1976), 271--272.

\bibitem{groth} \textrm{A.~Grothendieck}, 
\textit{Cohomologie locale des faisceaux coherents et th\'eor\`emes de Lefschetz 
locaux et globaux}. 
SGA 2, Soci\'et\'e Math\'ematique de France (2005).

\bibitem{hal} \textrm{M.~Halic}, 
\textit{Vector bundles on projective varieties whose restrictions 
to an ample subvariety split}. preprint http://arxiv.org/abs/1401.4732.

\bibitem{hal+taj} \textrm{M.~Halic, R.~Tajarod}, 
\textit{About the cohomological dimension of certain stratified varieties}. 
preprint. 

\bibitem{hart-ag} \textrm{R.~Hartshorne}, 
\textit{Algebraic Geometry}. $8^{\textrm{th}}$ ed. 
Springer-Verlag, New York, 1997. 

\bibitem{hart-as} \textrm{R.~Hartshorne}, 
\textit{Ample Subvarieties of Algebraic Varieties}. 
Lect. Notes Math. \textbf{156}, Springer-Verlag, Berlin, 1970.

\bibitem{ho} \textrm{G.~Horrocks},
\textit{Vector bundles on the punctured spectrum of a local ring}. 
Proc. Lond. Math. Soc., \textbf{14} (1964), 689--713.

\bibitem{kon} J.~Konarski, 
\textit{The B-B decomposition via Sumihiro's theorem}. 
J.~Algebra \textbf{182} (1996), 45--51.

\bibitem{ma-oe} \textrm{P.M.~Marques,  L.~Oeding}, 
\textit{Splitting criteria for vector bundles on the symplectic isotropic Grassmannian}. 
Matematiche (Catania) 64 (2009), 155--176. 

\bibitem{ma} \textrm{F.~Malaspina}, 
\textit{A few splitting criteria for vector bundles}. 
Ric. Mat. 57 (2008), no. 1, 55--64.

\bibitem{mat} \textrm{S.~Matsumura},
\textit{Asymptotic cohomology vanishing and a converse 
to the Andreotti-Grauert theorem on surfaces}. 
Ann. Inst. Fourier \textbf{63} (2013), 2199--2221.

\bibitem{mitch} \textrm{S.~Mitchell}, 
\textit{Quillen's theorem on buildings and the loops on a symmetric space}.
Enseign. Math. \textbf{34} (1988), 123--166.

\bibitem{ott} \textrm{G.~Ottaviani}, 
\textit{Some extensions of Horrocks criterion to vector bundles 
on Grassmannians and quadrics}. Ann. Mat. Pura Appl. 155 (1989), 317--341. 

\bibitem{ot} \textrm{J.~Ottem}, 
\textit{Ample subvarieties and $q$-ample divisors}. 
Adv. Math. \textbf{229} (2012), 2868--2887.

\bibitem{so} \textrm{A.~Sommese}, 
\textit{Submanifolds of abelian varieties}. 
Math. Ann. \textbf{233} (1978), 229--256.

\bibitem{sp} \textrm{T.A.~Springer}, 
\textit{Linear algebraic groups}, $2^{\rm nd}$ ed. 
Progress in Mathematics, vol. 9, Birkh\"auser, 1998. 

\bibitem{to} \textrm{B.~Totaro}, 
\textit{Line bundles with partially vanishing cohomology}. 
J. Eur. Math. Soc. \textbf{15} (2013), 731--754.

\end{thebibliography}
\end{document}